\RequirePackage{fix-cm}
\documentclass[smallextended]{svjour3}       % onecolumn (second format)
\smartqed  % flush right qed marks, e.g. at end of proof
\usepackage{graphicx,amssymb,graphicx,amsmath,xcolor,subfigure,caption}
\usepackage{epstopdf}
\usepackage{cite}
\usepackage{float}
\usepackage{tikz,setspace}
\usepackage{array}
\usepackage{hyperref}
\usepackage{booktabs}

\begin{document}

\title{Dynamics of an HIV/AIDS transmission model with protection awareness and fluctuations%\thanks{Grants or other notes
}


\author{Xuanpei Zhai \and
        Wenshuang Li \and
        Fengying Wei \and
        Xuerong Mao
}


\institute{Xuanpei Zhai \at
    School of Mathematics and Statistics,
    Fuzhou University, Fuzhou 350116, P.R. China \\
    \email{zhaixp2022@shanghaitech.edu.cn}
    \and
    Wenshuang Li \at
    School of Mathematics and Statistics,
    Fuzhou University, Fuzhou 350116, P.R. China \\
    \email{210320031@fzu.edu.cn}
    \and
    Fengying Wei \at
    School of Mathematics and Statistics,
    Fuzhou University, Fuzhou 350116, P.R. China \\
    Center for Applied Mathematics of Fujian Province,
    Fuzhou University, Fuzhou 350116, P.R. China  \\
    \email{weifengying@fzu.edu.cn}
    \and
    Xuerong Mao \at
    Department of Mathematics and Statistics,
    University of Strathclyde, Glasgow G1 1XH, UK \\
    \email{x.mao@strath.ac.uk}
}

\date{Received: date / Accepted: date}

\maketitle

\begin{abstract}
    We establish a stochastic HIV/AIDS model for the individuals with
    protection awareness and reveal how the protection awareness plays its important role in the control of AIDS. We
    firstly show that there exists a global positive solution for the
    stochastic model. By constructing Lyapunov functions, the ergodic
    stationary distribution when $R_{0}^{s}>1$ and the extinction when
    $R_{0}^{e}<1$ for the stochastic model are obtained. A number of numerical simulations by using positive preserving truncated Euler-Maruyama
    method (PPTEM) are performed to illustrate the theoretical results. Our new results show that the
    detailed publicity has great impact on the control of AIDS
    compared with the extensive publicity, while the continuous
    antiretroviral therapy (ART) is helpful in the control of HIV/AIDS.

\end{abstract}

\section{Introduction}
\label{intro} Infectious diseases caused one-quarter of the global
deaths \cite{ref1}. Various factors such as media campaigns,
population migration and temperature changes influenced the
spreading of infectious diseases. Since June 6 of 1981, the first
global case of HIV (Human Immunodeficiency Virus) infection was
announced, human beings have been fighting against HIV for over
four decades. As of the year 2022, about 37.7 million individuals
have been infected with HIV \cite{ref2}. As we have known today,
the transmission of HIV took place through blood, semen, cervical
(or, vaginal secretions) and breast milk as well. Especially, an
infected individual was unaware of the protection and was lack of
active treatment, HIV often broke down the immune system of the
infected individual and eventually turned into the acquired immune
deficiency syndrome (AIDS). Although there was no drug or vaccine
for HIV, the antiretroviral therapy (ART) could prolong the life
expectancy of an infected individual and make it
approach that of the uninfected individuals. Meanwhile, the
infected individuals with ART treatment do not retransmit HIV to
their sexual partners \cite{ref2}. From 2000 to 2018, the number
of new HIV infected individuals fell by 37\%, HIV-related deaths
fell by 45\%, and ART treatment saved 13.6 million individuals. At
the end of 2018, about 23.3 million individuals had received ART
treatment \cite{ref2}. Thus, UNAIDS puts forward the
90\%-90\%-90\% plan(90\% of AIDS infected people know they are
infected, 90\% of confirmed AIDS patients to be treated, and 90\%
of HIV in treated patients' body is suppressed), and set the great
goal of eliminating the AIDS epidemic by 2030 in UN General
Assembly Resolution 70/266.

Mathematical modelling of HIV/AIDS and its kinetic behavior
analysis can well predict the development trend of HIV/AIDS. Many
scholars have already studied the HIV/AIDS model and its kinetic
behavior. For example, Silva and Torres \cite{ref3} obtained the
results on the global stability of the HIV/AIDS model by
considering bilinear incidence rates. Ghosh et al. \cite{ref4}
studied the effect of medias and self-imposed psychological fears
on disease dynamics by separating the susceptible into the unaware
and the aware individuals. Later, Zhao et al. \cite{ref6} modified
the model established by Fatmawati et al. \cite{ref5}, and
considered piecewise fractional differential equations and
investigated the effect of protection awareness on HIV
transmission. However, the transmissions of infectious disease
were inevitably affected by the environmental noises in the real
circumstances. In other words, the numbers of the individuals in
each compartment were usually fluctuated due to the emergence of
infectious disease and control as well by local governments.
Therefore, the epidemic models with fluctuations in practice were
necessary to investigate their long-term dynamics. For instance,
Mao et al. \cite{ref7} found that the small fluctuations to the
deterministic models effectively suppress the rapid increment of
the population. Other recent epidemic models in \cite{ref8, ref9,
ref10, ref11,
    ref12, ref13, ref13.2, ref13.3} also governed the fluctuations to
describe the diversities of their models. More precisely,
\cite{ref13.2,ref13.3} found that small fluctuations produced the
long-term persistence, and large fluctuations led to the
extinction of infectious diseases in stochastic HIV/AIDS models.
Liu \cite{ref13.4} discovered that the higher order fluctuations
made HIV/AIDS eradicative under sufficient conditions. Meanwhile,
Wang \cite{ref13.5} also figured out the extinction and
persistence in the mean depended on the fluctuations of main
parameters.

We formulate a stochastic HIV/AIDS model with protection awareness
by considering the environmental noises into the model of
Fatmawati et al. \cite{ref5}. We next provide the expression of
the basic reproduction number for the deterministic model. In
Section 3, we prove theoretically that there exists a unique
positive solution to the stochastic model, and the existence of a
unique ergodic stationary distribution is investigated. Further,
we give a new threshold $R_0^e$ for the extinction of HIV/AIDS,
and the corresponding numerical simulations are demonstrated to
verify the theoretical results. Then we conclude that the detailed
publicity has great impact on the control of AIDS compared with
the extensive publicity. Moreover, the continuous antiretroviral
therapy (ART) is helpful to control the number of the individuals
with HIV/AIDS. Meanwhile, we figure out the main improvements
compared with other contributions, and give some suggestions to
control the long-term dynamics of HIV/AIDS.

\section{Establishment of the mathematical model}
\label{sec:2}

The number that an infected contacts with the susceptible per unit
of time is called the contact number, which is usually related to
the total population $N$ and denoted as $U(N)$. Let the
probability of infection per exposed individual be $\beta_0$, and
$\beta_0U(N)$ is called the effective exposure number, which
represents the infectivity of an infected individual per unit of
time. The total population is usually separated by the susceptible
individuals, the immune individuals, and the exposed individuals.
So, the proportion of the susceptible individuals in the whole
population is ${\frac{S(t)}{N(t)}}$, which is not being infected
by the infected individuals. Then, the number of the susceptible
being infected effectively at time $t$ is
$$
\beta_0U(N)\frac{S(t)}{N(t)}I(t)
$$
which is called the incidence rate. We assume in this paper that
the contact rate between the susceptible and the infected is
proportional to the total population, i.e., $U(N)=kN(t)$. Let
$\beta=\beta_0k$, then the incidence rate is rewritten as $\beta
S(t)I(t)$, which is called the bilinear incidence rate. Many
researchers have applied the bilinear incidence rates into their
HIV/AIDS models \cite{ref14,ref15,ref16,ref17} for further
discussions.

Fatmawati et al. \cite{ref5} developed a model of HIV/AIDS with
the protection awareness. Precisely, they divided the population
into five different groups, that is, the susceptible without
protection awareness $(S_u)$, the susceptible with protection
awareness $(S_a)$, the infected without ART treatment $(I)$, the
infected with ART treatment $(C)$ and the infected who eventually
developed into AIDS $(A)$. And, the total people size $N$ is
expressed as
$$
N=S_a+S_u+I+C+A.
$$
We simplified the model of Fatmawati et al. \cite{ref5} by considering bilinear incidence:
$$
\begin{array}{lll}
    \dot{S_{u}}(t)&=&\Lambda-\beta S_{u} I-(\alpha+\mu) S_{u}, \\
    \dot{S_{a}}(t)&=&\alpha S_{u}-(1-\varepsilon) \beta S_{a} I-\mu S_{a}, \\
    \dot{I}(t)&=&\beta I\left(S_{u}+(1-\varepsilon) S_{a}\right)+\eta C+v A-(\rho+\gamma+\mu) I, \\
    \dot{C}(t)&=&\rho I-(\eta+\mu) C, \\
    \dot{A}(t)&=&\gamma I-(v+\delta+\mu) A,
\end{array}\eqno{(2.1)}
$$
where $\Lambda$ is the recruitment rate, $\mu$ is the natural
death rate, and $\delta$ is the mortality rate of AIDS. $\alpha$
is the migration rate from $S_u$ to $S_a$, $\beta$ is the HIV
transmission rate, and $\varepsilon$ is the infection rate of
$S_a$. The parameters $\eta$, $\nu$, $\gamma$, and $\rho$
represent the transmission rates from $C$ to $I$, $A$ to $I$, $I$
to $A$ and $I$ to $C$. The mutual migrating mechanisms of each group
in the model could be demonstrated clearly in Figure 2.1. We suppose
that all parameters of model (2.1) are positive initiated with
$$
\begin{array}{ll}
    &S_u(0)=S_{u0}\geq0, S_a(0)=S_{a0}\geq0, I(0)=I_0\geq0,\\
    &C(0)=C_0\geq0, A(0)=A_0\geq0.
\end{array}
$$

\begin{figure}
    \centering
    \includegraphics[width=0.63\textwidth]{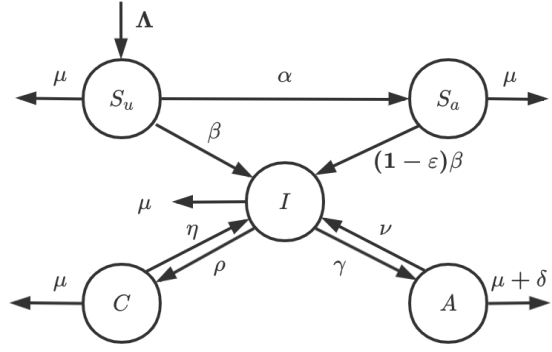}
    \caption*{Figure 2.1 Propagation mechanism}
\end{figure}

Now, we add five equations of model (2.1) and then we get
$$
\dot{N}(t)=\Lambda-\mu N-\delta A\leq \Lambda -\mu N.
$$
By comparing theorems, the positive invariant set of model (2.1)
is derived
$$
\Omega=\left\{\left(S_u,S_a,I,C,A\right)\in\mathbb R_+^5:0\le
N\le\frac{\Lambda}{\mu}\right\}.
$$
We only consider the biological properties of model (2.1) in the
set $\Omega$. The basic regeneration number of
model (2.1) can be obtained from the next generation matrix
approach \cite{ref17.1,ref17.2} as follows:
$$
\begin{array}{ll}
    R_0=\displaystyle\frac{\beta\left[\mu+\left(1-\varepsilon\right)\alpha\right]
        k_1k_2\Lambda}{\mu(\mu+\alpha)\left[\mu\left(k_2\left(k_1+\gamma\right)+\rho k_1+\gamma\delta\right)+\eta\gamma\delta\right]},
\end{array}
\eqno{(2.2)}
$$
with $k_1=\mu+\delta+\nu$, $k_2=\mu+\eta$.

Similarly to the proof of Theorem 2 in Fatmawati et al. \cite{ref5}, when $R_0<1$, we
obtain that model (2.1) has a unique boundary equilibrium point
$$
P_0=\left(\frac{\Lambda}{\mu+\alpha},\frac{\alpha\Lambda}{\mu\left(\mu+\alpha\right)},0,0,0\right),
$$
which is locally asymptotically stable in the set $\Omega$. We
also provide the expression of an endemic equilibrium point of
model (2.1) when $R_0>1$, that is
$$
P^\ast=(S_u^\ast,S_a^\ast,I^\ast,C^\ast,A^\ast),
$$
where
$$
S_u^\ast=\frac{\Lambda}{\beta I^\ast+\alpha+\mu},
S_a^\ast=\frac{\alpha\Lambda}{(\left(1-\varepsilon\right)\beta I^\ast+\mu)(\beta I^\ast+\alpha+\mu)},
$$
$$
C^\ast=\frac{\rho I^\ast}{k_2}, A^\ast=\frac{\gamma I^\ast}{k_1}.
$$
Substituting $S_u^\ast$, $S_a^\ast$, $C^\ast$, $A^\ast$ into the
third equation of (2.1) and making the left side be zero, then we
can get
$$
\begin{aligned}
    f\left(I^\ast\right)=&\beta I^\ast\left(\frac{\Lambda}{\beta
        I^\ast+\alpha+\mu}+
    \frac{(1-\varepsilon)\alpha\Lambda}{\left(\left(1-\varepsilon\right)\beta
        I^\ast+\mu\right)
        \left(\beta I^\ast+\alpha+\mu\right)}\right)+\eta\frac{\rho I^\ast}{k_2}\\
    &+\upsilon\frac{\gamma I^\ast}{k_1}-\left(\rho+\gamma+\mu\right)I^\ast\\
    =&I^\ast\left(\beta\Lambda
    k_1k_2\left(\left(1-\varepsilon\right)\beta
    I^\ast+\mu\right)+\beta\Lambda
    k_1k_2\left(1-\varepsilon\right)\alpha\right.\\&\left.
    +\left(k_1\eta\rho+k_2\upsilon\gamma-k_1k_2\left(\rho+\gamma+\mu\right)\right)\left(\beta
    I^\ast+\alpha+\mu\right)\left(\left(1-\varepsilon\right)\beta
    I^\ast+\mu\right)\right)\\&
    \left(k_1k_2\left(\beta I^\ast+\alpha+\mu\right)\left(\left(1-\varepsilon\right)\beta I^\ast+\mu\right)\right)^{-1}\\
    =&I^\ast (A^\ast{I^\ast}^2+B^\ast
    I^\ast+C^\ast)\left(k_1k_2\left(\beta
    I^\ast+\alpha+\mu\right)\left(\left(1-\varepsilon\right)\beta
    I^\ast+\mu\right)\right)^{-1},
\end{aligned}
$$
where
$$
\begin{aligned}
    P_1=&\left(1-\varepsilon\right)\beta^2\left(k_1\eta\rho+k_2\upsilon\gamma
    -k_1k_2\left(\rho+\gamma+\mu\right)\right)\\
    =&-\left(1-\varepsilon\right)\beta^2
    \left(\mu\upsilon\rho+\mu\delta\rho+\mu^2\rho+\eta\delta\rho+\eta\mu\rho+\mu\delta\gamma
    +\mu^2\gamma\right.\\
    &\left.+\mu\left(\nu+\delta+\mu\right)\left(\eta+\mu\right)\right)<0,\\
    P_2=&\beta\left(\left(1-\varepsilon\right)\alpha+\left(2-\varepsilon\right)\mu\right)
    \left(k_1\eta\rho+k_2\upsilon\gamma-k_1k_2\left(\rho+\gamma+\mu\right)\right)\\
    &+k_1k_2\Lambda\beta^2\left(1-\varepsilon\right),\\
    P_3=
    &k_1k_2\Lambda\beta\left(\mu+\left(1+\varepsilon\right)\alpha\right)
    +\mu\left(\mu+\alpha\right)\left(\ k_1\eta\rho+k_2\upsilon\gamma
    -k_1k_2\left(\rho+\gamma+\mu\right)\right)\\=
    &\mu\left(\mu+\alpha\right)\left(\mu\left[k_2\left(k_1+\gamma\right)+\rho k_1+\gamma\delta\right]
    +\eta\gamma\delta\right)\left(R_0-1\right).
\end{aligned}
$$
Obviously $f(0)=0$. When $R_0>1$, we have $P_1<0$ and $P_3>0$.
According to the Descartes sign rule, $f(I^\ast)$ has a unique
positive real root regardless of the sign of $P_2$. Therefore, the
endemic equilibrium point $P^\ast$ exists.

Motivated by the models described by stochastic differential
equations in \cite{ref18, ref19, ref20, ref21, ref22}, we
introduce the environmental fluctuations into model (2.1) similar
to that of Evans\cite{ref22.1} and Tan\cite{ref22.2}. We set
that environmental fluctuations are multiplicative white noise
types which are proportional to $S_u, S_a, I, C, A$. When $\Delta
t\rightarrow 0 $, we consider a Markov process 
$$
X(t) = \left(S_u(t),S_a(t), I(t), C(t), A(t)\right)^{T}
$$ 
with the following descriptions:
$$
\begin{aligned}
    &\mathbb{E}[S_u(t+\Delta t)-S_u(t)|X_t=x]\thickapprox \left[\Lambda-\beta S_{u}(t) I(t)-(\alpha+\mu) S_{u}(t)\right]\Delta t,\\
    &\mathbb{E}[S_a(t+\Delta t)-S_a(t)|X_t=x]\thickapprox\left[\alpha S_{u}(t)-(1-\varepsilon) \beta S_{a}(t) I(t)-\mu S_{a}(t)\right] \Delta t,\\
    &\mathbb{E}[I(t+\Delta t)-I(t)|X_t=x]\thickapprox \left[\beta I(t)\left(S_{u}(t)+(1-\varepsilon) S_{a}(t)\right)+\eta C(t)+vA(t) \right.\\ &\left. \hspace*{12.55em} -(\rho+\gamma+\mu) I(t) \right] \Delta t,\\
&\mathbb{E}[C(t+\Delta t)-C(t)|X_t=x]\thickapprox [\rho I(t)-(\eta+\mu) C(t)]\Delta t,\\
&\mathbb{E}[A(t+\Delta t)-A(t)|X_t=x]\thickapprox [\gamma
I(t)-(v+\delta+\mu) A(t)] \Delta t,
\end{aligned}
$$
and
$$
\begin{aligned}
&\mbox{Var}[S_u(t+\Delta t)-S_u(t)|X_t=x]\thickapprox\sigma_{1}^2 S_{u}^2(t)\Delta t,\\
&\mbox{Var}[S_a(t+\Delta t)-S_a(t)|X_t=x]\thickapprox\sigma_{2}^2 S_{a}^2(t)\Delta t,\\
&\mbox{Var}[I(t+\Delta t)-I(t)|X_t=x]\thickapprox\sigma_{3}^2 I^2(t)\Delta t,\\
&\mbox{Var}[C(t+\Delta t)-C(t)|X_t=x]\thickapprox\sigma_{4}^2 C^2(t)\Delta t,\\
&\mbox{Var}[A(t+\Delta t)-A(t)|X_t=x]\thickapprox\sigma_{5}^2 A^2(t)\Delta t.\\
\end{aligned}
$$
We therefore derive a stochastic epidemic model as follows:
$$
\begin{array}{lll}
    \mbox{d}S_{u}(t)&=&\left[\Lambda-\beta S_{u}(t) I(t)-(\alpha+\mu) S_{u}(t)\right] \mbox{d}t+\sigma_{1} S_{u}(t) \mbox{d}B_{1}(t), \\
    \mbox{d}S_{a}(t)&=&\left[\alpha S_{u}(t)-(1-\varepsilon) \beta S_{a}(t) I(t)-\mu S_{a}(t)\right] \mbox{d}t
    +\sigma_{2} S_{a}(t) \mbox{d}B_{2}(t), \\
    \mbox{d}I(t)&=&\left[\beta I(t)\left(S_{u}(t)+(1-\varepsilon) S_{a}(t)\right)+\eta C(t)+v A(t) \right.\\
    & &\left. -(\rho+\gamma+\mu) I(t) \right] \mbox{d}t +\sigma_{3} I(t) \mbox{d}B_{3}(t), \\
    \mbox{d}C(t)&=& [\rho I(t)-(\eta+\mu) C(t)] \mbox{d}t+\sigma_{4} C(t) \mbox{d}B_{4}(t), \\
    \mbox{d}A(t)&=& [\gamma I(t)-(v+\delta+\mu) A(t)] \mbox{d}t+\sigma_{5} A(t) \mbox{d}B_{5}(t).
\end{array}\eqno{(2.3)}
$$

Let $i=1,2,3,4,5$, then $B_i(t)$ are independent standard Brownian
motions with the initial values $B_i\left(0\right)=0$ and
$\sigma_i^2>0$ are the intensities of white noises, and the
initial values $ X(0)=(S_u(0), S_a(0), I(0), C(0), A(0))^{T}$ as
well.

\section{Existence and uniqueness of positive solution}
\label{sec:3} Let $X(t)$ be a homogeneous markov process in
$\mathbb R^d$, which satisfies the stochastic differential
equation
$$
\mbox{d}X(t)=b(x)\mbox{d}t+\sum_{r=1}^{k}{g_r\left(x\right)\mbox{d}B_r(t)},
\quad k\leq d.
$$
The diffusion matrix is defined as
$$
A(X)=(a_{ij}(X))_{d\times d}, \quad
a_{ij}(X)=\sum_{r=1}^{k}{g_r^i(X)g_r^j(X).}
$$

\vskip5pt \noindent\textbf{Lemma 3.1} \cite{ref23} Assume that
there exists a bounded open region $G\subset \mathbb R^d$ with
regular boundaries $\Gamma$, and it has the following properties:

(i) The minimum eigenvalue of the diffusion matrix $A(X)$ is
non-zero in its domain $G$ and one of its neighborhoods.

(ii) The average time $\tau$ for the path from $z$ to set $G$ is
finite when $z\in \mathbb R^d\backslash G$, and $\sup_{z\in
K}E^z\tau<\infty$ holds for each compact subset $K\in \mathbb
R^d$.

Then, the Markov process $X(t)$ has a unique ergodic stationary
distribution $\pi(\cdot )$. Let $f(X)$ be an integrable function
of $\pi$, for all $X\in \mathbb R^d$, the following formula holds:
$$
\mathbb{P}\left\{\lim_{t\rightarrow\infty}{\frac{1}{t}\int_{0}^{t}f\left(X(s)\right)\mbox{d}s=\int_{\mathbb
R^d}{f(X)\pi(\mbox{d}X)}}\right\}=1.
$$

\vskip5pt \noindent\textbf{Remark 3.1} The proof of Lemma 3.1 can
be found on pages 106-109 of Khasminskii \cite{ref23}. If there
exists a positive $M$ such that
$$
    \sum_{i,j=1}^{d}{a_{ij}\left(X\right)\xi_i\xi_j}\geq M\left|\xi\right|^2, \xi\in \mathbb R^d,
$$
then property (i) holds.

\vskip5pt Now, we will give two useful results, Lemma 3.2 and
Lemma 3.3, by using of Theorem 2.1 and Theorem 3.1 in
\cite{ref27}. In fact, we write down the conclusions without
consider the details of the proofs.

\vskip5pt \noindent\textbf{Lemma 3.2} \textit{ Let $X(t)$ be a
solution of $(2.3)$ initiated with $X_0\in \mathbb{R}_+^5$, then
$$
    \lim_{t\rightarrow\infty}{\frac{1}{t}\left(S_u(t)+S_a(t)+I(t)+C(t)+A(t)\right)}=0,
$$
and
$$
\begin{aligned}
    &\lim_{t\rightarrow\infty}{\frac{S_u(t)}{t}}=0 ,\lim_{t\rightarrow\infty}{\frac{S_a(t)}{t}}=0,
    \lim_{t\rightarrow\infty}{\frac{I(t)}{t}}=0,\\&\lim_{t\rightarrow\infty}{\frac{C(t)}{t}}=0,
    \lim_{t\rightarrow\infty}{\frac{A(t)}{t}}=0 \hspace{1em} \text{a.s.}
\end{aligned}
$$
}

\vskip5pt \noindent\textbf{Lemma 3.3} \textit{ Suppose that
$\mu>(\sigma_1^2\vee \sigma_2^2\vee \sigma_3^2 \vee \sigma_4^2\vee
\sigma_5^2)/2$, let $X(t)$ be a solution of $(2.3)$ initiated with
$X_0\in \mathbb{R}_+^5$. Then
$$
\begin{aligned}
    &\lim_{t\rightarrow\infty}{\frac{1}{t}\int_{0}^{t}{S_u(s)}\mbox{d}B_1(s)}=0,
    \lim_{t\rightarrow\infty}{\frac{1}{t}\int_{0}^{t}{S_a(s)}\mbox{d}B_2(s)}=0,\\
    &\lim_{t\rightarrow\infty}{\frac{1}{t}\int_{0}^{t}{I(s)}\mbox{d}B_3(s)}=0,
    \lim_{t\rightarrow\infty}{\frac{1}{t}\int_{0}^{t}{C(s)}\mbox{d}B_4(s)}=0,\\
    &\lim_{t\rightarrow\infty}{\frac{1}{t}\int_{0}^{t}{A(s)}\mbox{d}B_5(s)}=0 \hspace{1em} \text{a.s.}
\end{aligned}
$$
}

Before we start to study the dynamical behaviors of the stochastic
epidemic model (2.3), the existence of a global positive solution
is of importance. Next, we show that there exists a unique global
positive solution to (2.3) for any given initial value.

\vskip5pt \noindent\textbf{Theorem 3.1} \textit{ Model $(2.3)$ has
a unique global positive solution $X(t)\in \mathbb{R}_+^5$
initiated with $X_0\in \mathbb{R}_+^5$ for any $t\geq 0$. }

\vskip5pt \noindent\textbf{\emph{Proof}} It is obvious to check
that the local Lipschitz condition is satisfied for model (2.3)
initiated with $X_0\in \mathbb{R}_+^5$, so there exists a unique
local solution $X(t)$ for $t\in[0, \tau_e)$. To prove that $X(t)$
is global, our work is to verify $\tau_e=\infty$. Indeed, let
$n_0>1$ be large enough satisfying each component of $X(t)$ lies
in $[{\frac{1}{n_0}}, n_0]$. Define the stopping time
$$
\begin{array}{lll}
    {\tau_n}=&\inf \left\{t\in[0,\tau_e):\min{\left\{S_u(t),S_a(t),I(t),C(t),A(t)\right\}}\le\displaystyle\frac{1}{n} \right.\\
    &\quad \left.\text{or} \hspace{1mm} \max{\left\{S_u(t),S_a(t),I(t),C(t),A(t)\right\}}\geq n \right\}
\end{array}
$$
for any integer $n>n_0$. Let $\inf\varnothing=\infty$. As
$n\rightarrow\infty$, it is obvious that
$\left\{\tau_n\right\}_{n\geq n_0}$ is monotonically increasing.
We set $\tau_\infty=\lim_{n \to \infty}{\tau_n}$, then we get
$\tau_\infty\le\tau_e$ by the definition of stopping time. We
claim that $\tau_\infty=\infty$. What we claim is checked, which
ends the proof. By contradiction, there exists a pair of positvie
constants $T>0$ and $\varepsilon\in(0, 1)$ such that the
probability that $\tau_\infty\le T$ is larger than $\varepsilon$.
We rewrite as $\mathbb{P}\left\{\tau_n\le
T\right\}\geq\varepsilon$ for $n\geq  n_0$. Define a
$C^2$-function $U:\mathbb{R}_+^5\rightarrow \mathbb{R}_+$ by
$$
    \begin{aligned}
        U\left(x\right)=&\left(S_u-\theta-\theta \ln\frac{S_u}{\theta}\right)
        +\left(S_a-1-\ln S_a\right)+\left(I-1-\ln I\right)\\& +\left(C-1-\ln C\right)+\left(A-1-\ln A\right),
    \end{aligned}
$$
where $\theta\in \mathbb{R}_+$. By the scalar It\^{o}'s formula,
we get
$$
    \begin{aligned}
        \mbox{d}U\left(x\right)=&\mathcal{L}U\left(x\right)\mbox{d}t
        +\sigma_1\left(S_u-\theta\right)\mbox{d}B_1(t)+\sigma_2\left(S_a-1\right)\mbox{d}B_2(t)\\
        &+\sigma_3\left(I-1\right)\mbox{d}B_3(t)+\sigma_4\left(C-1\right)\mbox{d}B_4(t)+\sigma_5\left(A-1\right)\mbox{d}B_5(t),
    \end{aligned}
$$
where $\mathcal{L}U\left(x\right):\mathbb{R}_+^5\rightarrow \mathbb{R}_+$ is
$$
    \begin{aligned}
        \mathcal{L}U\left(x\right)= & - \frac{\theta\Lambda}{S_u}+\theta\beta I+\theta\left(\alpha+\mu\right)
        +\frac{\theta}{2}{\sigma_1^2}-\frac{\alpha S_u}{S_a}+\left(1-\varepsilon\right)\beta I\\
        &+\mu + \frac{1}{2}\sigma_2^2-\beta\left(S_u+\left(1-\varepsilon\right)S_a\right)-\frac{\eta C}{I}-\frac{\nu A}{I}\\
        &+\rho+\gamma+\mu
        +\frac{1}{2}\ \sigma_3^2-\frac{\rho I}{C}+\eta+\mu+\frac{1}{2}\ \sigma_4^2-\frac{\gamma I}{A}\\
        &+\nu+\delta+\mu+\frac{1}{2}\sigma_5^2+\Lambda\\
        &-\mu\left(S_u+S_a+I+C+A\right)-\delta A\notag
    \end{aligned}
$$
$$
\begin{aligned}
    \le &\ (\theta+1-\varepsilon)\beta I-\mu I+\theta\left(\alpha+\mu+\frac{1}{2}\sigma_1^2\right)+\Lambda+4\mu\\
    &+\rho+\gamma+\eta+\nu+\delta+\frac{1}{2}\left(\sigma_2^2{+\sigma}_3^2{+\sigma}_4^2{+\sigma}_5^2\right).\notag
\end{aligned}
$$
We let $\theta=\frac{\mu}{\beta}-1+\varepsilon$, then
$$
    \begin{aligned}
        \mathcal{L}U\left(x\right)\le
        &\left(\frac{\mu}{\beta}-1+\varepsilon\right)
        \left(\alpha+\mu+\frac{1}{2}\sigma_1^2\right)+\Lambda+4\mu+\rho+\gamma+\eta+\nu+\delta\\
        &+\frac{1}{2}\left(\sigma_2^2{+\sigma}_3^2{+\sigma}_4^2{+\sigma}_5^2\right):=Q,\\
    \end{aligned}
$$
which further gives
$$
    \begin{aligned}
        \mbox{d}U\left(x\right)\le&Q\mbox{d}t+\sigma_1\left(S_u-\theta\right)\mbox{d}B_1(t)
        +\sigma_2\left(S_a-1\right)\mbox{d}B_2(t)+\sigma_3\left(I-1\right)\mbox{d}B_3(t)\\
        &+\sigma_4\left(C-1\right)\mbox{d}B_4(t)+\sigma_5\left(A-1\right)\mbox{d}B_5(t).
    \end{aligned}\eqno{(3.1)}
$$
Integrating (3.1) from $0$ to $\tau_n\land T=\min{\left\{\tau_n,
T\right\}}$, taking the expectation, we get
$$
     \mathbb{E} \left(U\left(x(\tau_n\land  T)\right)\right)
     \le U\left(X_0\right)+Q\mathbb{E}\left(\tau_n\land  T\right)\le U\left(X_0\right)+QT.\eqno{(3.2)}
$$
When $n\geq n_0$, let $\Omega_n=\{\tau_n\le T\}$, then the
inequality $\mathbb{P}\{{\tau_n\le T}\}\geq\varepsilon$ transforms
into $\mathbb{P}\{{\Omega_n}\}\geq\varepsilon$. For each
$\omega\in\Omega_n$, $S_u$ takes value $n-\theta-\theta
\ln\frac{n}{\theta}$ or $\frac{1}{n}-\theta+\theta\ln{\left(\theta
n\right)}$ at time $\tau_n\land T$, so do $S_a$, $I$, $C$ and $A$.

Obviously, inequality (3.2) can be transformed into
$$
    \begin{aligned}
        U(X_0)+QT&\geq \mathbb{E} \left[\mathbb{I}_{\Omega_n(\omega)}U\left(x\left(\tau_n\land  T\right)\right)\right]\\
        &\geq \varepsilon\left[\left(n-\theta-\theta \ln\frac{n}{\theta}\ \right)
        \land \left(\frac{1}{n}-\theta+\theta\ln{\left(\theta n\right)}\right)\right],\notag
    \end{aligned}
$$
where $\mathbb{I}_{\Omega_n\left(\omega\right)}$ is the index
function of $\Omega_n\left(\omega\right)$. Let $n\rightarrow
\infty$, we get
$$
    \infty>U\left(X_0\right)+QT=\infty.
$$
This is a contradiction. The proof is complete.

\section{Existence of a unique ergodic stationary distribution}
\label{sec:4} Sufficient conditions for the existence of
stationary distribution and ergodicity of model (2.3) are given
below, which also implies that HIV/AIDS is persistent in the mean.
We demote the stochastic index
$$
    R_0^s=\frac{\beta\left[\mu+\left(1-\varepsilon\right)\alpha\right]
    k_1k_2\Lambda}{k_5k_6\left(k_1k_2\left(k_4+\frac{1}{2}\sigma_3^2\right)
    -\left(k_1+\frac{1}{2}\sigma_5^2\right)\rho\eta-\left(k_2+\frac{1}{2}\sigma_4^2\right)\gamma\nu\right)},\notag
$$
with
$$
\begin{aligned}
    & k_1=\mu+\delta+\nu, k_2=\mu+\eta,  k_3={\frac{\beta\left(\mu+\left(1-\varepsilon\right)\alpha\right)\Lambda}{\mu\left(\mu+\alpha\right)}},\\
    & k_4=\rho+\gamma+\mu, k_5=\mu+\frac{1}{2}\sigma_2^2,
    k_6=\mu+\alpha+\frac{1}{2}\sigma_1^2.
\end{aligned}
$$
When $\sigma_1=\sigma_2=\sigma_3=\sigma_4=\sigma_5=0$, $R_0^s$
degenerates to $R_0$ in (2.2).

\vskip5pt \noindent\textbf{Theorem 4.1}\textit{ Model $(2.3)$ has
a unique stationary distribution, and it is ergodic when
$R_0^s>1$.}

\vskip5pt \noindent\textbf{\emph{Proof}} The diffusion matrix
$$
D\left(x\right)=\text{diag}
\left\{\sigma_1^2S_u^2,\sigma_2^2S_a^2{,\sigma}_3^2I^2{,\sigma}_4^2C^2,\sigma_5^2A^2\right\}
$$
of model (2.3) is positive definite, then condition (i) is clearly
established. Therefore, we only need to prove that condition (ii)
holds. First, we create a $C^2$-function
$G:\mathbb{R}_+^5\rightarrow \mathbb{R}_+$ by
$$
    \begin{aligned}
        G\left(x\right)=&M (-c_1\ln S_u-c_2\ln S_a+c_3k_1k_2\beta\frac{1}{\alpha}S_a-c_3k_1k_2\ln I \\
        &+c_3\rho\eta \ln A+c_3\gamma\nu \ln C)+\left(S_u+S_a+I+C+A\right)^{m+1}\\
        &-\ln S_u-\ln S_a-\ln C-\ln A\\
        := & MV_1+V_2+V_3+V_4+V_5+V_6,
    \end{aligned}
    \notag
$$
where $c_i\in\mathbb{R}_+$ for $i=1,2,3$; $M>0$ is a sufficiently
large positive number, $m$ is a sufficiently small positive
number, $M$ and $m$ satisfy
$$
\begin{aligned}
     &-3M\left(\sqrt[3]{R_0^s}-1\right)+(M\beta^\ast+\left(2-\varepsilon\right)\beta)\varepsilon_1+B+\alpha+\mu\\
     &+\frac{1}{2}\sigma_1^2+\mu+\frac{1}{2}\sigma_2^2+k_2+\frac{1}{2}\sigma_4^2+k_1+\frac{1}{2}\sigma_5^2\le-2,
\end{aligned}\eqno{(4.1)}
$$
and
$$
    \mu-\frac{1}{2}m\left(\sigma_1^2\lor \sigma_2^2\lor \sigma_3^2\lor \sigma_4^2\lor \sigma_5^2\right)>0,\eqno{(4.2)}
$$
here $\beta^\ast$ and $B$ are defined in (4.3) and (4.4)
respectively.

We obtained that
$$
    \lim_{n\rightarrow \infty}\ \inf_{\left(S_u,S_a,I,C,A\right)\in \mathbb{R}_+^5\backslash U} G\left(S_u,S_a,I,C,A\right)=+\infty,
$$
where
$$
    U=\left(\frac{1}{n},n\right)\times\left(\frac{1}{n},n\right)\times\left(\frac{1}{n},n\right)
    \times\left(\frac{1}{n},n\right)\times\left(\frac{1}{n},n\right).
$$
Since $G\left(S_u,S_a,I,C,A\right)$ is a continuous function,
there must be a minimum value $\widetilde{G}$. Define a
non-negative $C^2$-function
$$
    V_1\left(S_u,S_a,I,C,A\right)=G\left(S_u,S_a,I,C,A\right)-\widetilde{G},
$$
then we apply the It\^{o}'s formula on $V_1$:
$$
    \begin{aligned}
        \mathcal{L}V_1=
        &-c_1\left(\frac{\Lambda}{S_u}-\beta I-\left(\alpha+\mu\right)\right)+\frac{1}{2}c_1\sigma_1^2
        -c_2\left(\frac{\alpha S_u}{S_a}-\left(1-\varepsilon\right)\beta I-\mu\right)\\
        &+\frac{1}{2}c_2\sigma_2^2-c_3k_1k_2\left(\beta\left({S_u}+\left(1-\varepsilon\right)S_a\right)
        +\frac{\eta C}{I}+\frac{\nu A}{I}-k_4\right)\\
        &+\frac{1}{2}c_3k_1k_2\sigma_3^2+c_3\rho\eta\left(\frac{\gamma I}{A}-k_1\right)
        -\frac{1}{2}c_3\rho\eta\sigma_5^2+c_3\gamma\nu\left(\frac{\rho I}{C}-k_2\right)\\
        &-\frac{1}{2}c_3\gamma\nu\sigma_4^2+c_3k_1k_2\beta\frac{1}{\alpha}\left(\alpha S_u-\left(1-\varepsilon\right)\beta S_aI-\mu S_a\right)\\
        =&-\left(c_1\frac{\Lambda}{S_u}+c_2\frac{\alpha S_u}{S_a}+c_3k_1k_2\beta\left(\frac{\mu}{\alpha}+1-\varepsilon\right)S_a\right)
        +c_1\left(\alpha+\mu+\frac{1}{2}\sigma_1^2\right)\\
        &+c_2\left(\mu+\frac{1}{2}\sigma_2^2\right)+c_3\left[k_1k_2\left(k_4+\frac{1}{2}\sigma_3^2\right)\
         -\left(k_1+\frac{1}{2}\sigma_5^2\right)\rho\eta \right.\\
        &\left.-\left(k_2+\frac{1}{2}\sigma_4^2\right)\gamma\nu\right]+\Big(c_1\beta+c_2\left(1-\varepsilon\right.)\beta\left.
        +c_3\rho\gamma\left(\frac{\eta}{A}+\frac{\nu}{C}\right)\right)I\\
        &-c_3k_1k_2\left(\frac{\eta C}{I}+\frac{\nu A}{I}+\left(1-\varepsilon\right)\beta^2{\frac{1}{\alpha}S_a}I\right).\notag
    \end{aligned}
$$
Using $a+b+c\geq3\sqrt[3]{abc}$ for positive $a,b$ and $c$, we get
$$
    \begin{aligned}
        \mathcal{L}V_1&\le-3\sqrt[3]{c_1c_2c_3\beta\left[\mu+\left(1-\varepsilon\right)\alpha\right]k_1k_2\Lambda}
        +c_1\left(\alpha+\mu+\frac{1}{2}\sigma_1^2\right)\\&+c_2\left(\mu+\frac{1}{2}\sigma_2^2\right)
        +c_3\left(k_1k_2\left(k_4+\frac{1}{2}\sigma_3^2\right)\ -\left(k_1+\frac{1}{2}\sigma_5^2\right)\rho\eta\right.\\&
        \left.-\left(k_2+\frac{1}{2}\sigma_4^2\right)\gamma\nu\right)+\left(c_1\beta+c_2\left(1-\varepsilon\right)\beta
        +c_3\rho\gamma\left(\eta+\nu\right)\right)I,\notag
    \end{aligned}
$$
we let
$$
    c_1=\frac{1}{\alpha+\mu+\frac{1}{2}\sigma_1^2},\quad
    c_2=\frac{1}{\mu+\frac{1}{2}\sigma_2^2},\\
$$
$$
    c_3=\frac{1}{k_1k_2\left(k_4+\frac{1}{2}\sigma_3^2\right) -\left(k_1+\frac{1}{2}\sigma_5^2\right)\rho\eta
    -\left(k_2+\frac{1}{2}\sigma_4^2\right)\gamma\nu},
$$
$$
    \beta^\ast=c_1\beta+c_2\left(1-\varepsilon\right)\beta+c_3\rho\gamma\left(\eta+\nu\right),
$$
$$
    \hat M=\mu-\frac{1}{2}m\left(\sigma_1^2\vee \sigma_2^2\vee \sigma_3^2\vee \sigma_4^2\vee \sigma_5^2\right),
$$
and then
$$
    \mathcal{L}V_1\le-3\left(\sqrt[3]{R_0^s}-1\right)+\beta^\ast I. \eqno{(4.3)}
$$
Similarly,
$$
    \begin{aligned}
        \mathcal{L}V_2=&\left(m+1\right) N^m\left(\Lambda-\mu N-\delta A\right)\\
        &+\frac{1}{2}\left(m+1\right)m N^{m-1}(\sigma_1^2S_u^2+\sigma_2^2S_a^2+\sigma_3^2I^2+\sigma_4^2C^2+\sigma_5^2A^2)\\
        \le&\left(m+1\right) N^m\left(\Lambda-\mu N\right)\\
        &+\frac{1}{2}\left(m+1\right)m N^{m+1}(\sigma_1^2\vee \sigma_2^2\vee \sigma_3^2\vee \sigma_4^2\vee \sigma_5^2)\\
        \le&\left(m+1\right)\Lambda N^m-\left(m+1\right)\hat MN^{m+1}\\
        \le&B-\frac{1}{2}\left(m+1\right)\hat M\left(S_u^{m+1}+S_a^{m+1}+I^{m+1}+C^{m+1}+A^{m+1}\right),
    \end{aligned}
    \eqno{(4.4)}
$$
where
$$
    B=\sup_{N\in\left(0,\infty\right)}{{}\left(m+1\right)\left(\Lambda N^m-\frac{1}{2}\hat MN^{m+1}\right)}<\infty.
$$
We thus derive
\begin{align}\nonumber
    \mathcal{L}V_3&=-\frac{1}{S_u}\left[\Lambda-\beta S_uI-\left(\alpha+\mu\right)S_u\right]+\frac{1}{2S_u^2}(\sigma_1^2S_u^2)\\
    &=-\frac{\Lambda}{S_u}+\beta I+\left(\alpha+\mu\right)+\frac{1}{2}\sigma_1^2, \tag{4.5}\\
    \mathcal{L}V_4&=-\frac{\alpha S_u}{S_a}+\left(1-\varepsilon\right)\beta I+\mu+\frac{1}{2}\sigma_2^2, \tag{4.6}\\
    \mathcal{L}V_5&=-\frac{\rho I}{C}+k_2+\frac{1}{2}\sigma_4^2, \tag{4.7}\\
    \mathcal{L}V_6&=-\frac{\gamma I}{A}+k_1+\frac{1}{2}\sigma_5^2. \tag{4.8}
\end{align}
From (4.3)-(4.8) we can get
$$
    \begin{aligned}
        \mathcal{L}V\le&-3M\left(\sqrt[3]{R_0^s}-1\right)+\left(M\beta^\ast+\left(2-\varepsilon\right)\beta\right)I
        +B-\frac{\Lambda}{S_u}+\alpha+\mu\\
        &-\frac{1}{2}\left(m+1\right)\hat M(S_u^{m+1}+S_a^{m+1}+I^{m+1}+C^{m+1}+A^{m+1})\\
        &+\frac{1}{2}\sigma_1^2-\frac{\alpha S_u}{S_a}+\mu+\frac{1}{2}\sigma_2^2-\frac{\rho I}{C}
        +k_2+\frac{1}{2}\sigma_4^2-\frac{\gamma I}{A}+k_1+\frac{1}{2}\sigma_5^2.
    \end{aligned}
    \eqno{(4.9)}
$$
We define a bounded region
$$
\begin{aligned}
    H=&\left\{X(t)\in \mathbb{R}_+^5:\varepsilon_1\le S_u(t)\le\frac{1}{\varepsilon_1},
    \varepsilon_1^2\le S_a(t)\le\frac{1}{\varepsilon_1^2},\varepsilon_1\le I(t)\le\frac{1}{\varepsilon_1}, \right.\\
    &\left.\varepsilon_1^2\le C(t) \le\frac{1}{\varepsilon_1^2}{,\varepsilon}_1^2\le A(t)\le\frac{1}{\varepsilon_1^2}\right\},
\end{aligned}
$$
where $\varepsilon_1>0$ is sufficiently small and satisfies:
\begin{align}
    &-\frac{\Lambda}{\varepsilon_1}+F\le-1,\tag{4.10}\\
    &-\frac{\alpha}{\varepsilon_1}+F\le-1,\tag{4.11}\\
    &-3M\left(\sqrt[3]{R_0^s}-1\right)+(M\beta^\ast+\left(2-\varepsilon\right)\beta)\varepsilon_1+F-e\le-1,\tag{4.12}\\
    &-\frac{\rho}{\varepsilon_1}+F\le-1,\tag{4.13}\\
    &-\frac{\gamma}{\varepsilon_1}+F\le-1,\tag{4.14}\\
    &-\frac{1}{2}\left(m+1\right)\hat M\frac{1}{\varepsilon_1^{m+1}}+F\le-1,\tag{4.15}\\
    &-\frac{1}{2}\left(m+1\right)\hat M\frac{1}{\varepsilon_1^{2m+2}}+F\le-1,\tag{4.16}\\
    &-\frac{1}{4}\left(m+1\right)\hat M\frac{1}{\varepsilon_1^{m+1}}+F\le-1,\tag{4.17}
\end{align}
combined with (4.1), we denote
$$
    \begin{aligned}
        F&:=B+e+\alpha+\mu+\frac{1}{2}\sigma_1^2+\mu+\frac{1}{2}\sigma_2^2+k_2+\frac{1}{2}\sigma_4^2+k_1+\frac{1}{2}\sigma_5^2,\\
        e&=\sup_{I\in(0,\infty)}{\left\{-\frac{1}{4}\left(m+1\right)\hat M I^{m+1}
        +\left(M\beta^\ast+\left(2-\varepsilon\right)\beta\right)I\right\}}<\infty.
    \end{aligned}
$$
Obviously $\mathbb{R}_+^5\backslash H=D_1\cup D_2\cup \cdots\cup
D_{10}$, where
$$
     \begin{aligned}
        D_1&=\left\{X(t)\in \mathbb{R}_+^5:0<S_u<\varepsilon_1\right\},\\
        D_2&=\left\{X(t)\in \mathbb{R}_+^5:0<S_a<\varepsilon_1^2,S_u\geq\varepsilon_1\right\},\\
        D_3&=\left\{X(t)\in \mathbb{R}_+^5:0<I<\varepsilon_1\right\},\\
        D_4&=\left\{X(t)\in \mathbb{R}_+^5:0<C<\varepsilon_1^2,I\geq\varepsilon_1\right\},\\
        D_5&=\left\{X(t)\in \mathbb{R}_+^5:0<A<\varepsilon_1^2,I\geq\varepsilon_1\right\},\\
        D_6&=\left\{X(t)\in \mathbb{R}_+^5:S_u\geq{1/\varepsilon}_1\right\},\\
        D_7&=\left\{X(t)\in \mathbb{R}_+^5:S_a\geq1/\varepsilon_1^2\right\},\\
        D_8&=\left\{X(t)\in \mathbb{R}_+^5:I\geq{1/\varepsilon}_1\right\},\\
        D_9&=\left\{X(t)\in \mathbb{R}_+^5:C\geq1/\varepsilon_1^2\right\},\\
        D_{10}&=\left\{X(t)\in \mathbb{R}_+^5:A\geq1/\varepsilon_1^2\right\}.
    \end{aligned}
$$
We next discuss each case as follows:

\noindent\textbf{Case 1}. When $X(t)\in D_1$, according to (4.1),
(4.9), (4.10), we can get
$$
    \mathcal{L}V\le-\frac{\Lambda}{S_u}+F\le-\frac{\Lambda}{\varepsilon_1}+F\le-1.
$$
\noindent\textbf{Case 2}. When $X(t)\in D_2$, according to (4.1),
(4.9), (4.11), we can get
$$
    \mathcal{L}V\le-\frac{\alpha S_u}{S_a}+F\le-\frac{\alpha}{\varepsilon_1}+F\le-1.
$$
\noindent\textbf{Case 3}. When $X(t)\in D_3$, according to (4.1),
(4.9), (4.12), we can get
$$
    \begin{aligned}
        \mathcal{L}V\le&-3M\left(\sqrt[3]{R_0^s}-1\right)+[M\beta^\ast+\left(2-\varepsilon\right)\beta]I+F-e\\
        \le&-3M\left(\sqrt[3]{R_0^s}-1\right)+[M\beta^\ast+\left(2-\varepsilon\right)\beta]\varepsilon_1+F-e\\\le&-1.
    \end{aligned}
$$
\noindent\textbf{Case 4}. When $X(t)\in D_4$, according to (4.1),
(4.9), (4.13), we can get
$$
    \mathcal{L}V\le-\frac{\rho I}{C}+F\le-\frac{\rho}{\varepsilon_1}+F\le-1.
$$
\noindent\textbf{Case 5}. When $X(t)\in D_5$, according to (4.1),
(4.9), (4.14), we can get
$$
    \mathcal{L}V\le-\frac{\gamma I}{A}+F\le-\frac{\gamma}{\varepsilon_1}+F\le-1.
$$
\noindent\textbf{Case 6}. When $X(t)\in D_6$, according to (4.1),
(4.9), (4.15), we can get
$$
    \begin{aligned}
        \mathcal{L}V\le&-\frac{1}{2}\left(m+1\right)\hat MS_u^{m+1}+F\\
        \le&-\frac{1}{2}\left(m+1\right)\hat M\frac{1}{\varepsilon_1^{m+1}}+F\le-1.
    \end{aligned}
$$
\noindent\textbf{Case 7}. When $X(t)\in D_7$, according to (4.1),
(4.9), (4.16), we can get
$$
    \begin{aligned}
        \mathcal{L}V\le&-\frac{1}{2}\left(m+1\right)\hat MS_a^{m+1}+F\\
        \le&-\frac{1}{2}\left(m+1\right)\hat M\frac{1}{\varepsilon_1^{2m+2}}+F\le-1.
    \end{aligned}
$$
\noindent\textbf{Case 8}. When $X(t)\in D_8$, according to (4.1),
(4.9), (4.17), we can get
$$
    \begin{aligned}
        \mathcal{L}V\le&-\frac{1}{2}\left(m+1\right)\hat MI^{m+1}+(\hat M\beta^\ast+\left(2-\varepsilon\right)\beta)I+F-e\\
        \leq&-\frac{1}{4}\left(m+1\right)\hat MI^{m+1}+F\\
        \le&-\frac{1}{4}\left(m+1\right)\hat M\frac{1}{\varepsilon_1^{m+1}}+F\le-1.
    \end{aligned}
$$
\noindent\textbf{Case 9}. When $X(t)\in D_9$, according to (4.1),
(4.9), (4.16), we can get
$$
    \begin{aligned}
        \mathcal{L}V\le&-\frac{1}{2}\left(m+1\right)\hat MC^{m+1}+F\\
        \le&-\frac{1}{2}\left(m+1\right)\hat M\frac{1}{\varepsilon_1^{2m+2}}+F\le-1.
    \end{aligned}
$$
\noindent\textbf{Case 10}. When $X(t)\in D_{10}$, according to
(4.1), (4.9), (4.16), we can get
$$
    \begin{aligned}
        \mathcal{L}V\le&-\frac{1}{2}\left(m+1\right)\hat MA^{m+1}+F\\
        \le&-\frac{1}{2}\left(m+1\right)\hat M\frac{1}{\varepsilon_1^{2m+2}}+F\le-1.
    \end{aligned}
$$
Therefore, we get $\mathcal{L}V\le-1$ as $X(t)\in
\mathbb{R}_+^5\backslash H$. So condition (ii) of Lemma 3.1 is
performed. The proof is complete.

\section{Extinction}
\label{sec:8} In this section, we will establish the sufficient
conditions for the extinction of infectious disease HIV/AIDS.
Denote
$$
    \left\langle S_u(t)\right\rangle=\frac{1}{t}\int_{0}^{t}{S_u(s)}\text ds,
    \ \left\langle S_a(t)\right\rangle=\frac{1}{t}\int_{0}^{t}{S_a(s)}\text ds.
$$

\vskip5pt \noindent\textbf{Theorem 5.1}\textit{ Suppose that
$\mu>(\sigma_1^2\vee\sigma_2^2\vee\sigma_3^2\vee\sigma_4^2\vee\sigma_5^2)/2$,
for any initial value $X_0\in \mathbb{R}_+^5$, if
$$
    R_0^e=\frac{\beta\left[\mu+\left(1-\varepsilon\right)\alpha\right]\Lambda}{\mu+\frac{1}{3}\hat \sigma}<1,
$$
with
$$
    \hat \sigma = \left(\delta+\frac{1}{2}\sigma_5^2\right)\wedge \frac{1}{2}\sigma_3^2\wedge
    \frac{1}{2}\sigma_4^2,
$$
then HIV/AIDS will become extinct, and the solution of model
$(2.3)$ satisfies:
$$
    \lim_{t\rightarrow\infty}{\frac{1}{t}\ln{\left(I(t)+C(t)+A(t)\right)}}<0.
$$
Moreover,
$$
    \begin{aligned}
         &\lim_{t\rightarrow\infty}\left\langle S_u(t)\right\rangle=\frac{\Lambda}{\mu+a},
         \ \lim_{t\rightarrow\infty}\left\langle S_a(t)\right\rangle=\frac{\alpha\Lambda}{\mu(\mu+a)}.\\
         &\lim_{t\rightarrow\infty}I(t)=0,\ \lim_{t\rightarrow\infty}C(t)=0,\
         \lim_{t\rightarrow\infty}A(t)=0.
    \end{aligned}
$$
}

\vskip5pt \noindent\textbf{\emph{Proof}} It is easy to check that
\begin{align}
    \text{d}S_u(t)&\le\left[\Lambda-\left(\alpha+\mu\right)S_u(t)\right]\text{d}t+\sigma_1S_u(t)\mbox{d}B_1(t),\tag{5.1}\\
    \text{d}S_a(t)&\le\left[\alpha S_u(t)-\mu S_a(t)\right]\text{d}t+\sigma_2S_a(t)\mbox{d}B_2(t).\tag{5.2}
\end{align}
By Lemma 3.2 and Lemma 3.3, after integration, we have
\begin{align}
     \lim_{t\rightarrow\infty}\frac{1}{t}\left(S_u(t)-S_u\left(0\right)\right)
     \le\lim_{t\rightarrow\infty}\left[\Lambda-\left(\alpha+\mu\right)\left\langle S_u(t)\right\rangle
     +\frac{\sigma_1}{t}\int_{0}^{t}{S_u(s)}\mbox{d}B_1(s)\right],\notag
\end{align}
which further shows
\begin{align}
\lim_{t\rightarrow\infty}\left\langle S_u(t)\right\rangle \le
\frac{\Lambda}{\mu+a}.\tag{5.3}
\end{align}
By the similar argument, we get
\begin{align}
\lim_{t\rightarrow\infty}\frac{1}{t}\left(S_a(t)-S_a\left(0\right)\right)
\le\lim_{t\rightarrow\infty}\left[\alpha\left\langle
S_u(t)\right\rangle-\mu\left\langle S_a(t)\right\rangle
+\frac{\sigma_2}{t}\int_{0}^{t}{S_a(s)}\mbox{d}B_2(s)\right],\notag
\end{align}
which thus implies
\begin{align}
    \lim_{t\rightarrow\infty}\left\langle S_a(t)\right\rangle\le\frac{\alpha\Lambda}{\mu(\mu+a)}.\tag{5.4}
\end{align}
Here we define $P\left(X\right)=I+C+A$, so It\^{o}'s formula gives
that
\begin{align}
    \mbox{d}\ln P\left(X\right)=
    w\mbox{d}t\notag
    +\frac{1}{P\left(X\right)}\left[\sigma_3I\mbox{d}B_3(t)+\sigma_4C\mbox{d}B_4(t)+\sigma_5A\mbox{d}B_5(t)\right],\tag{5.5}
\end{align}
with
$$
\begin{aligned}
    w=&\frac{1}{P\left(X\right)}\left(\beta I\left(S_u+\left(1-\varepsilon\right)S_a\right)-\mu P\left(X\right)-\delta A\right)\notag\\
    &-\frac{1}{2P^2\left(X\right)}\left(\sigma_3^2I^2+\sigma_4^2C^2+\sigma_5^2A^2\right).
\end{aligned}
$$
Due to the facts that
\begin{equation*}
     -\frac{\delta A}{P\left(X\right)}\le-\frac{\delta A^2}{P^2\left(X\right)},\
     \frac{I}{P\left(X\right)}\le1,\frac{C}{P\left(X\right)}\le1,\frac{A}{P\left(X\right)}\le1,
\end{equation*}
then (5.5) is simplified as
\begin{align}
    w\le&\beta\left(S_u+\left(1-\varepsilon\right)S_a\right)-\mu-\frac{\delta A^2}{P^2\left(X\right)}
    -\frac{1}{2P^2\left(X\right)}\left(\sigma_3^2I^2+\sigma_4^2C^2+\sigma_5^2A^2\right)\notag\\
    =&\beta\left(S_u+\left(1-\varepsilon\right)S_a\right)-\mu-\frac{1}{P^2\left(X\right)}
    \left(\frac{1}{2}\sigma_3^2I^2+\frac{1}{2}\sigma_4^2C^2+(\delta+\frac{1}{2}\sigma_5^2)A^2\right)\notag \\
    \le&\beta\left(S_u+\left(1-\varepsilon\right)S_a\right)-\mu-\frac{I^2+C^2+A^2}{P^2\left(X\right)}\hat \sigma \notag\\
    \le&\beta\left(S_u+\left(1-\varepsilon\right)S_a\right)-\mu-\frac{1}{3}\hat \sigma. \tag{5.6}
\end{align}
The integration on (5.6) gives that
\begin{align}
    \frac{1}{t}\left(\ln P\left(X\right)-\ln P\left(X(0)\right)\right)\le
    &\beta\left\langle S_u(t)\right\rangle+\left(1-\varepsilon\right)\left\langle S_a(t)\right\rangle-\mu-\frac{1}{3}\hat \sigma \notag\\
    &+\frac{\sigma_3}{t}B_3(t)+\frac{\sigma_4}{t}B_4(t)+\frac{\sigma_5}{t}B_5(t).\tag{5.7}
\end{align}
Applying the strong law of numbers, we get
\begin{equation}
    \lim_{t\rightarrow\infty}{\frac{B_i(t)}{t}=0\quad\left(i=3,4,5\right).}\tag{5.8}
\end{equation}
When $t\rightarrow\infty$, by (5.3), (5.4), (5.8) and $R_0^e<1$,
(5.6) can be simplified as
\begin{equation*}
    \begin{aligned}
        \lim_{t\rightarrow\infty}\sup{\frac{\ln P\left(X\right)}{t}}&\le\beta\frac{\Lambda}{\mu+a}
        +\left(1-\varepsilon\right)\frac{\alpha\Lambda}{\mu\left(\mu+a\right)}-\mu\ -\frac{1}{3}\hat \sigma \\
        &=\left(R_0^e-1\right)\left(\mu+\frac{1}{3}\hat \sigma \right)<0.
    \end{aligned}
\end{equation*}
Therefore, we derive
\begin{equation}
     \lim_{t\rightarrow\infty}I(t)=0,\quad
     \lim_{t\rightarrow\infty}C(t)=0,\quad
     \lim_{t\rightarrow\infty}A(t)=0.\tag{5.9}
\end{equation}
Furthermore, we consider
\begin{equation*}
        N=\Lambda-\mu\left(N\right)-\delta A+\sigma_1S_u+\sigma_2S_a+\sigma_3I+\sigma_4C+\sigma_5A,
\end{equation*}
the integration implies that
\begin{equation*}
     \begin{aligned}
        &\frac{N(t)-N(0)}{t}\\
        &=\Lambda-\mu\left(\left\langle S_u(t)\right\rangle+\left\langle S_a(t)\right\rangle
        +\left\langle I(t)\right\rangle+\left\langle C(t)\right\rangle+\left\langle A(t)\right\rangle\right)
        -\delta\left\langle A(t)\right\rangle\\
        &\quad +\frac{\sigma_1}{t}\int_{0}^{t}{S_u(s)}\mbox{d}B_1(s)
        +\frac{\sigma_2}{t}\int_{0}^{t}{S_a(s)}\mbox{d}B_2(s)
        +\frac{\sigma_3}{t}\int_{0}^{t}I(s)\mbox{d}B_3(s)\\
        &\quad +\frac{\sigma_4}{t}\int_{0}^{t}C(s)\mbox{d}B_4(s)+\frac{\sigma_5}{t}\int_{0}^{t}A(s)\mbox{d}B_5(s).
    \end{aligned}
\end{equation*}
By Lemma 3.2 and Lemma 3.3, together with (5.3), (5.4) and (5.9),
the following expressions are obtained:
\begin{equation*}
     \lim_{t\rightarrow\infty}\left\langle S_u(t)\right\rangle=\frac{\Lambda}{\mu+a},\
     \lim_{t\rightarrow\infty}\left\langle S_a(t)\right\rangle=\frac{\alpha\Lambda}{\mu(\mu+a)}.
\end{equation*}
Thus, the proof of Theorem 5.1 is complete.

%%%%%%%%%%%%

\section{Examples and numerical simulations}
\label{sec:6}

We take the parameters in this section from Fatmawati et al.
\cite{ref5} except for the values for $\beta$ and $\sigma$. We
further govern the positive preserving truncated Euler-Maruyama
method (also referred as PPTEM) in \cite{ref28} to simulate the
long-term properties of the solution. Let
$$X(t_k) = (S_u(t_k), S_a(t_k), I(t_k), C(t_k), A(t_k)), k = 0,1,2,\cdots$$
be the discrete solution of model (2.3) with $t_k=k\Delta$, then the corresponding
discretization equations are written as
\begin{equation}
    \begin{split}
        S_u(t_{k+1}) =& S_u(t_k) + (\Lambda + f_{11} + f_{12})\Delta + g_{1} \sqrt{\Delta}, \\
        S_a(t_{k+1}) =& S_a(t_k) + (f_{21} + f_{22})\Delta + g_{2}\sqrt{\Delta}, \\
        I(t_{k+1}) =& I(t_k) + (f_{31} + f_{32})\Delta + g_{3}\sqrt{\Delta}, \\
        C(t_{k+1}) =& C(t_k) + f_{41}\Delta + g_{4} \sqrt{\Delta}, \\
        A(t_{k+1}) =& A(t_k) + f_{51}\Delta + g_{5} \sqrt{\Delta}, \\
    \end{split}
    \tag{6.1}
\end{equation}
for $k = 0, 1, 2, 3, \cdots$, and
\begin{equation}
    \begin{split}
        f_{11} =& - (\alpha + \mu) \hat\pi_0(S_u(t_k)),
        f_{12} = - \beta \hat\pi_0(S_u(t_k) I(t_k)), \\
        f_{21} =& \alpha \hat\pi_0(S_u(t_k)) - \mu \hat\pi_0(S_a(t_k)),\\
        f_{22} =& - (1-\varepsilon) \beta \hat\pi_0(S_a(t_k) I(t_k)),\\
        f_{31} =& \eta \hat\pi_0(C(t_k)) + \nu \hat\pi_0(A(t_k)) - (\rho + \gamma + \mu)\hat\pi_0(I(t_k)),\\
        f_{32} =& \beta \hat\pi_0(I(t_k) [\hat\pi_0(S_u(t_k)) + (1-\varepsilon)\hat\pi_0(S_a(t_k))],\\
        f_{41} =& \rho \hat\pi_0(I(t_k)),
        f_{51} = \gamma \hat\pi_0(I(t_k)),\\
        g_{1} =& \sigma r_{1,k} \hat\pi_0(S_u(t_k)), \hspace{3mm} g_{2} = \sigma r_{2,k} \hat\pi_0(S_a(t_k)),
        \hspace{3mm} g_{3} = \sigma r_{3,k} \hat\pi_0(I(t_k)),\\
        g_{4} =& \sigma r_{4,k} \hat\pi_0(C(t_k))), \hspace{3mm} g_{5} = \sigma r_{5,k} \hat\pi_0(A(t_k)),\\
    \end{split}
    \tag{6.2}
\end{equation}
where $\Delta$ is the stepsize and $r_{i,k} ~(i=1,2,3,4,5\text{ and } k=0,1,2,3,\cdots)$ are
independent random variables with the normal distribution
$\mathcal{N}(0, 1)$ and the function $\hat\pi_0(u)$ is defined as
$\hat\pi_0(u)=0 \vee u$. The discretization equations can be denoted as:
$$
X(t_{k+1}) = X(t_k) + [f_1 (X(t_k)) + f_2 (X(t_k))] \Delta +
g(X(t_k)) \Delta B_k
$$
where
\begin{equation}\nonumber
    \begin{split}
        &f_1 = (\Lambda + f_{11}, f_{21}, f_{31}, f_{41}, f_{51})^T, f_2 = (f_{12}, f_{22}, f_{32}, 0, 0)^T,\\
        &g = (g_1, g_2, g_3, g_4, g_5)^T.
    \end{split}
\end{equation}

We define a strictly increasing function $z:
\mathbb{R}_+\to\mathbb{R}_+$ by $z(u)=u$ for $u \geq 1$, which
gives the inverse function of $z^{-1}: [1, \infty)\to\mathbb{R}_+$
with the form $z^{-1}(u) = u$ for $u \geq 1$. We also define a
strictly decreasing function $h: (0, 1] \to [1, \infty)$ by
$h(\Delta) = \hat h \Delta^{-\frac{1}{4}}$ with $\hat h = 1 \vee
z(1) \vee |x(0)|$ and $|X(0)| = \sqrt{S_u^2(0) + S_a^2(0) + I^2(0)
    + C^2(0) + A^2(0)}$.

Define
\begin{equation}\nonumber
    \pi_\Delta(X) = \left(|X| \wedge z^{-1}(h(\Delta)) \right)\displaystyle\frac{X}{|X|},
\end{equation}
and
\begin{equation}\nonumber
    \hat \pi_\Delta(X) = (\Delta \vee S_u,\Delta \vee S_a,\Delta \vee I,\Delta \vee C,\Delta \vee A)^T.
\end{equation}
Let $\bar X_\Delta (t_k)$ be the intermediate step in order to get
nonnegative preserving truncated EM (NPTEM) solution $X_\Delta
(t_k)$, and $\bar X_\Delta (0) = X_\Delta (0) = X(0)$ be the
initial value, then the discretization equation of NPTEM is then
defined by
\begin{equation}
\bar X_\Delta (t_{k+1}) = \bar X_\Delta (t_k) + [f_1(\bar X_\Delta
(t_k)) + f_2(X_\Delta (t_k))] \Delta + g(\bar X_\Delta (t_k))
\Delta B_k, \tag{6.3}
\end{equation}
\begin{equation}
X_\Delta (t_{k+1}) = \hat \pi_0(\pi_\Delta(\bar X(t_{k+1}))),
\tag{6.4}
\end{equation}
for $k = 0, 1, 2, \cdots$, where $\Delta B_k = B(t_{k+1}) -
B(t_k)$, and then we extend the definition of $X_\Delta (\cdot)$
from the grid points $t_k$ to the whole $t \geq 0$ by defining
\begin{equation*}
x_\Delta (t) = x_\Delta (t_k) \hspace{1em} \text{for} \hspace{1em} t \in [t_k, t_{k+1}), \hspace{1mm} k=0,1,2,\cdots. \tag{6.5}
\end{equation*}
Together with (6.3) and (6.5), the positivity preserving truncated
EM (PPTEM) solution is consequently derived by
$$
X_\Delta^+ (t_{k+1}) = \hat \pi_\Delta(\pi_\Delta(\bar
X(t_{k+1}))), \hspace{1mm} k=0,1,2,\cdots.
$$

Next, we present the simulations in Indonesia and China by using
of PPTEM and predict the development and prevalence of HIV/AIDS
for next five decades.

\noindent\textbf{Example 6.1}

We firstly study the epidemics of HIV/AIDS  in Indonesia. We
choose $S_u(0)=129789089$, $S_a(0)=100000000$, $I(0)=7195$,
$C(0)=0$, $A(0)=3716$, $\Delta = 10^{-2}$ and let
$$\Lambda=\frac{229800000}{67.39}, \beta=\frac{0.3465}{229800000}, \mu=\frac{1}{67.39}, $$
and other parameters be
$$
\begin{aligned}
  &\alpha=0.2351, \varepsilon=0.3243, \eta=0.2059, \upsilon=0.7661, \\
  &\gamma=0.1882, \rho=0.00036523, \delta=0.7012.
\end{aligned}
$$

\noindent By Theorem 4.1, the stochastic index is $R_0^s \approx
2.075>1$ as $\sigma_i(i=1,2,3,4,5)= 0.05$, so HIV/AIDS is
persistent in a long run (see the left of Figure 6.1). The
population size in each compartment of the stochastic model (2.3)
fluctuates around the endemic equilibrium point
$$P^*=(12267874,85638867, 18584806, 30748, 2359917).$$
Furthermore, the solution of
model (2.3) has a unique stationary distribution, which is ergodic
when $R_0^s \approx 2.2676>1$ for $\sigma_i= 0.01~ (i=1,2,3,4,5)$
and $T=40000, 60000, 80000$ respectively, the population size in
each compartment is presented on the right of Figure 6.1.

\newpage
\begin{figure}[H]
    \centering
    \includegraphics[width=2in]{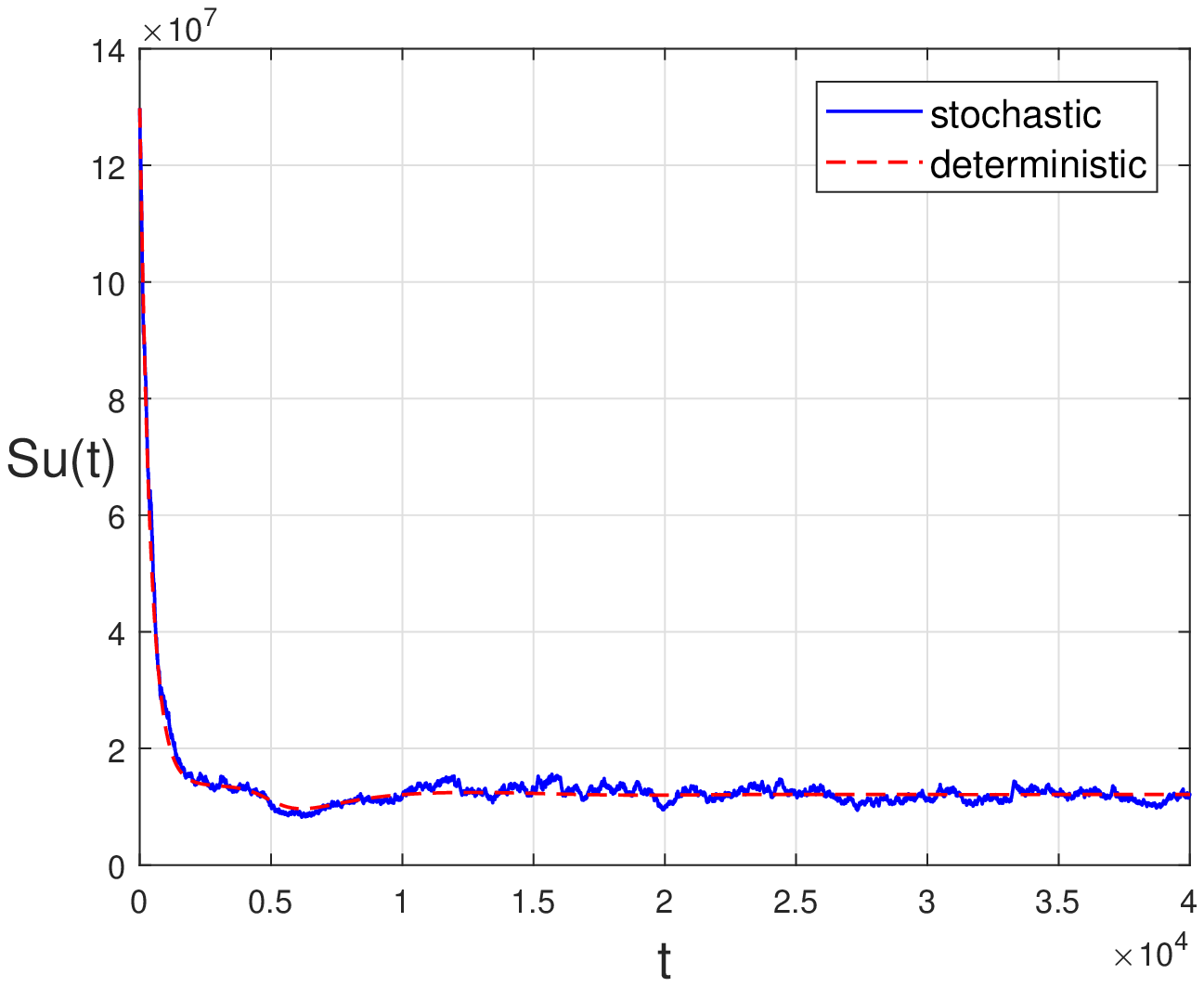}   \includegraphics[width=2in]{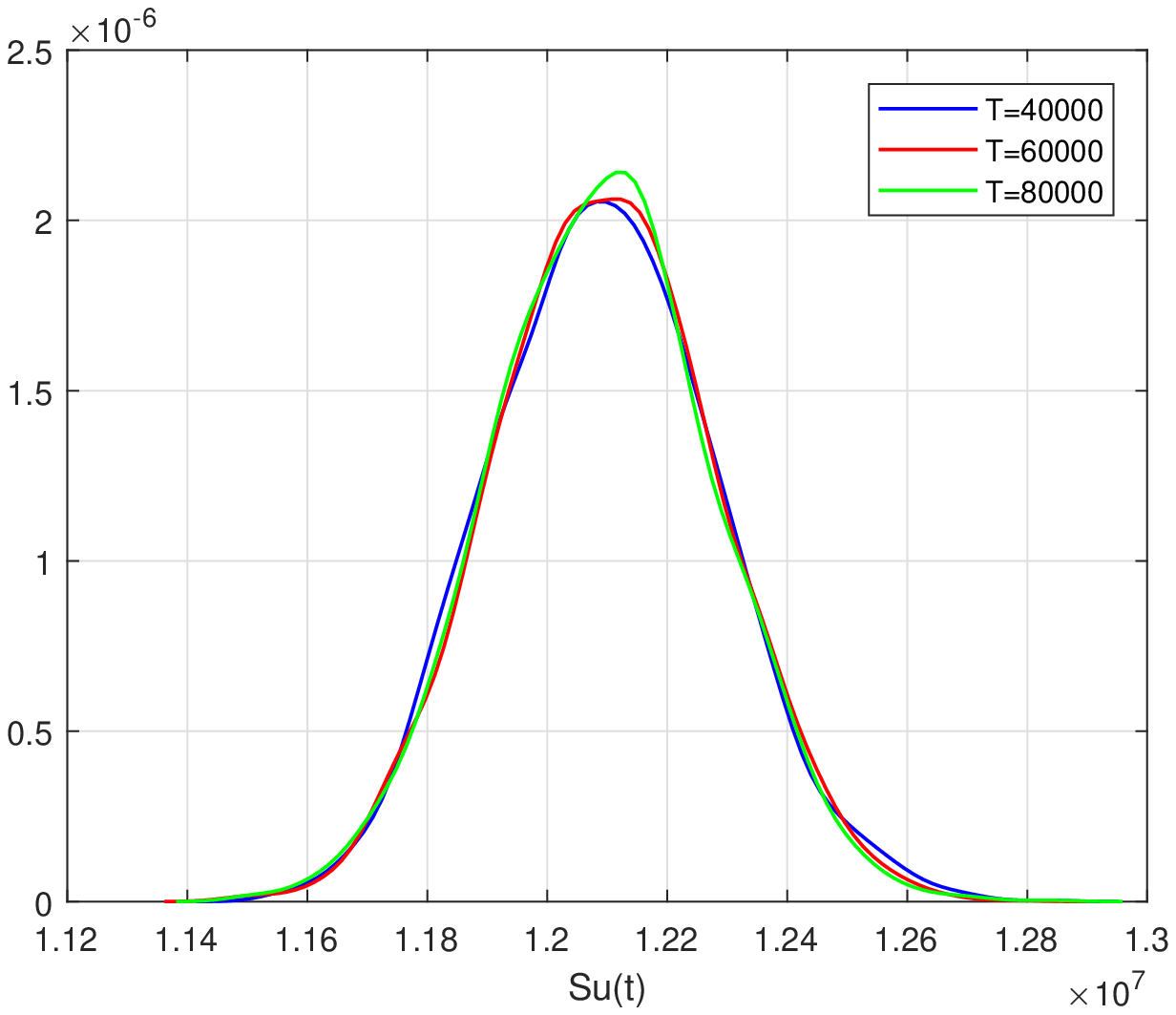}
    \includegraphics[width=2in]{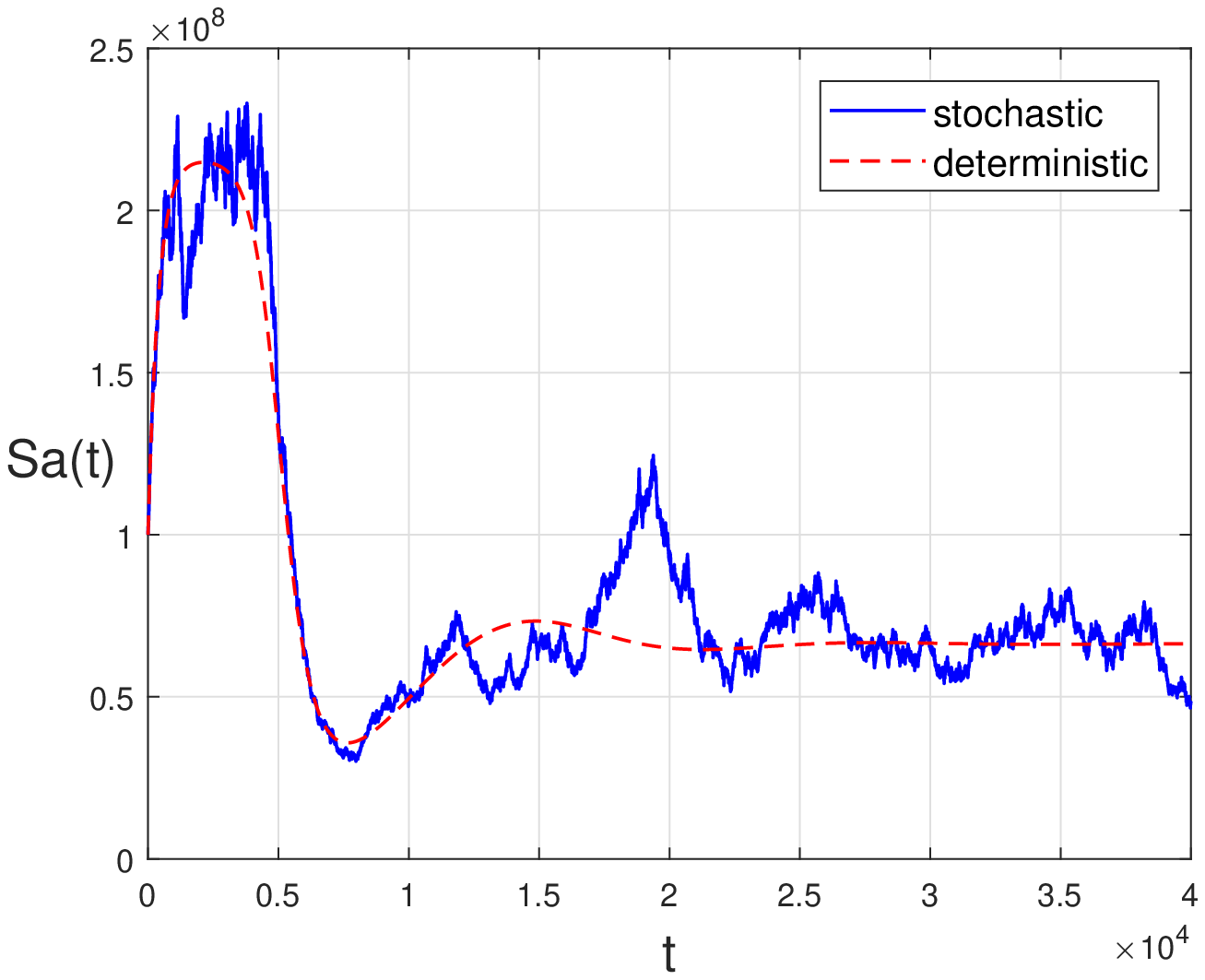}   \includegraphics[width=2in]{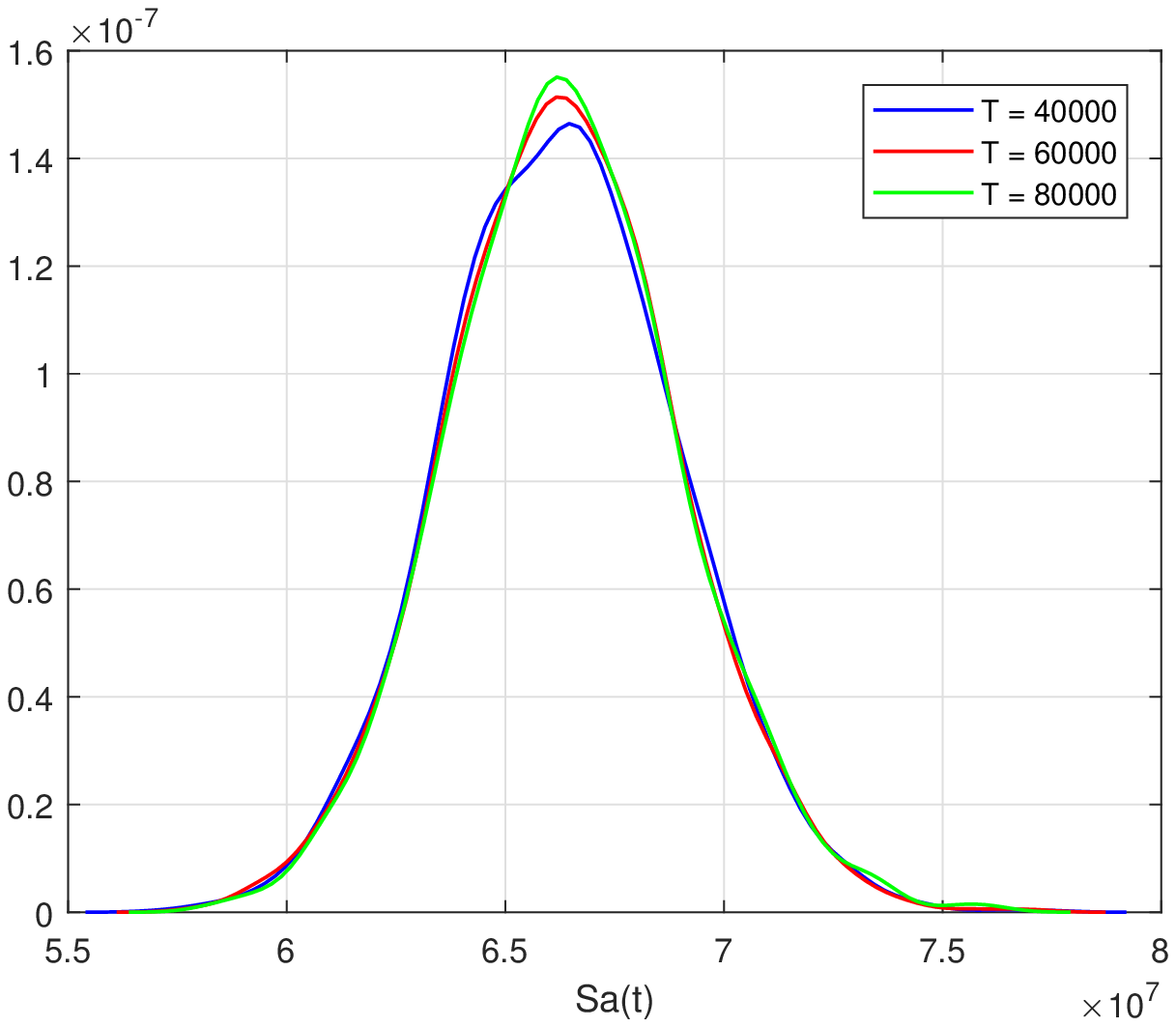}
    \includegraphics[width=2in]{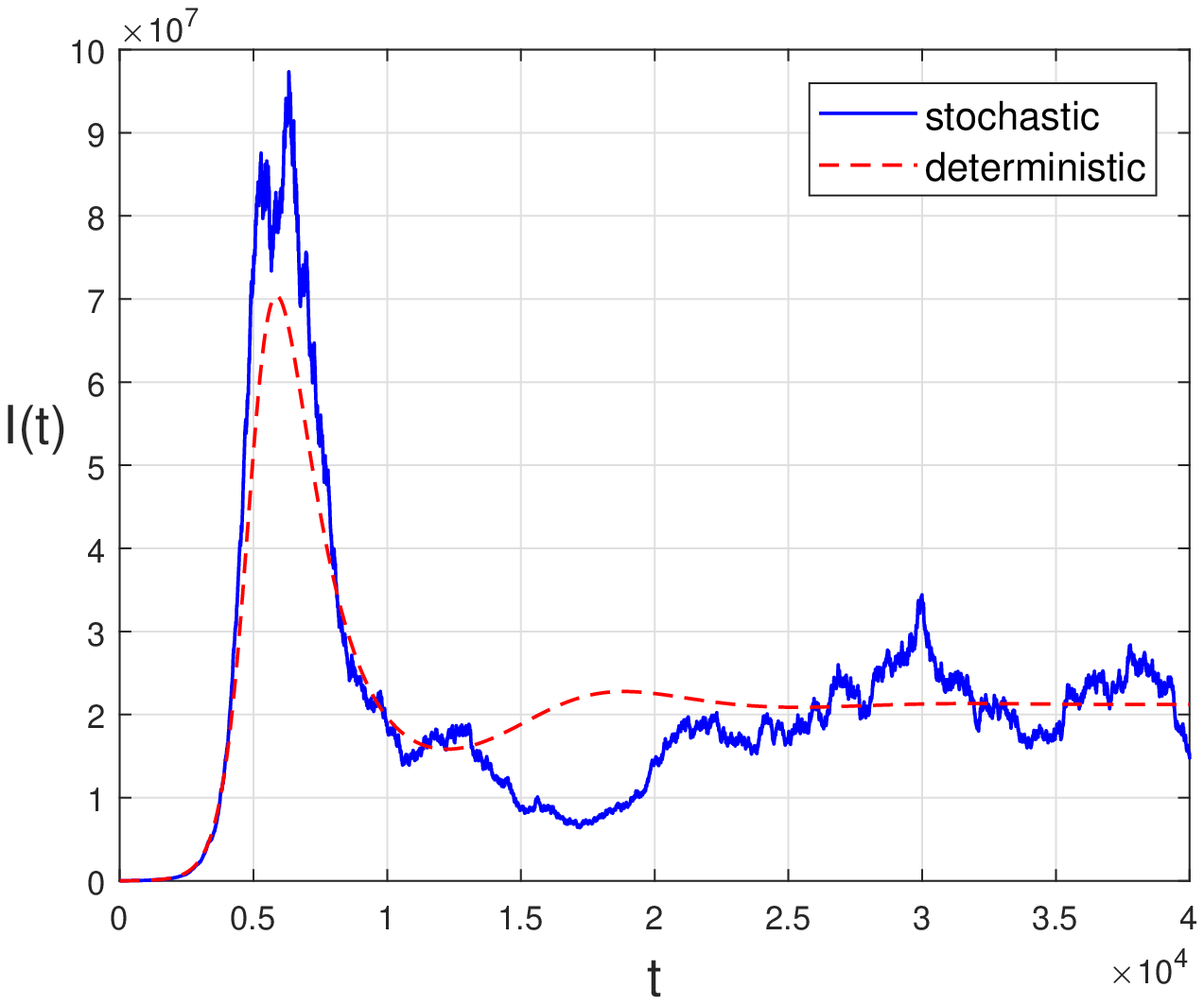}    \includegraphics[width=2in]{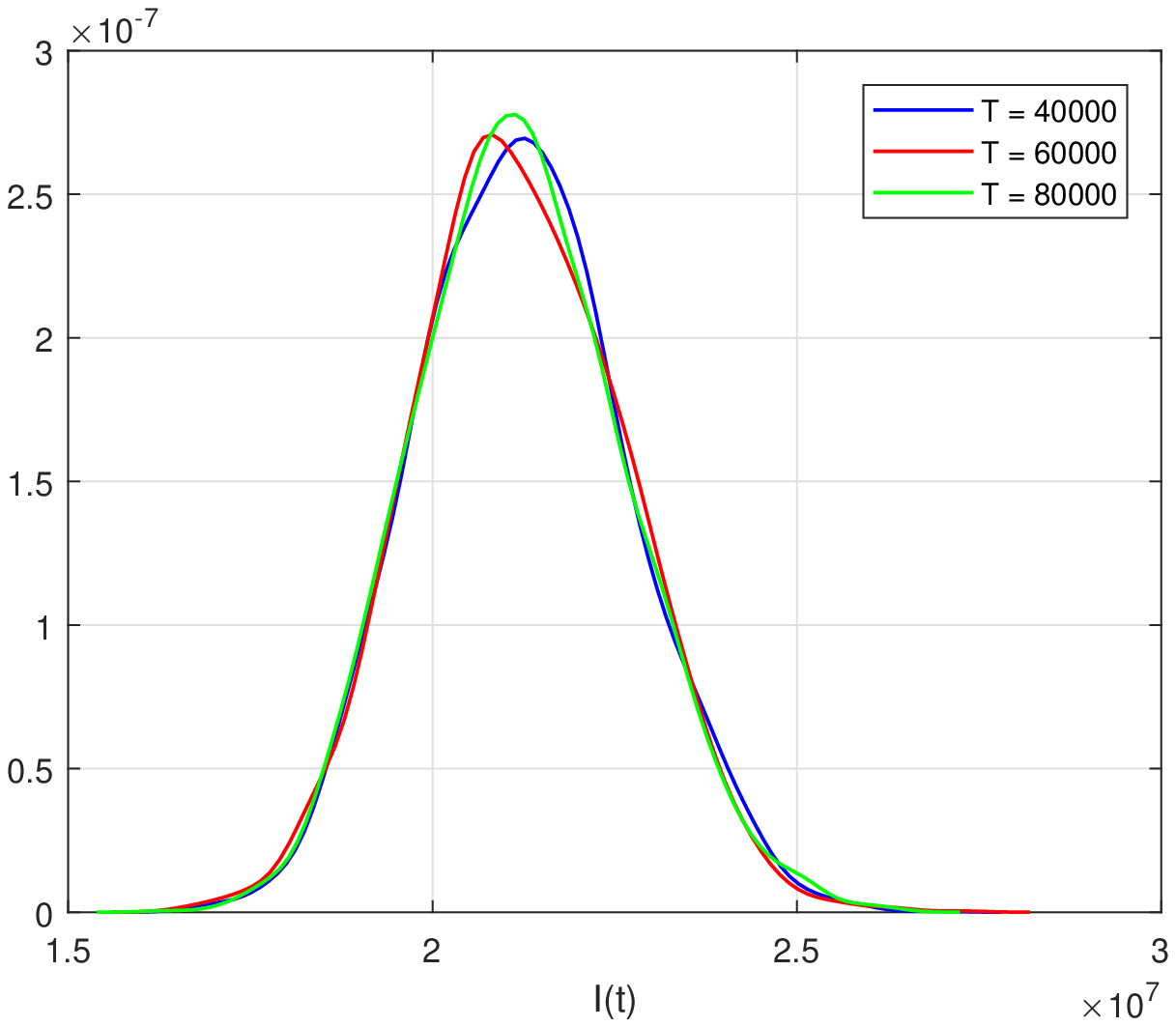}
    \includegraphics[width=2in]{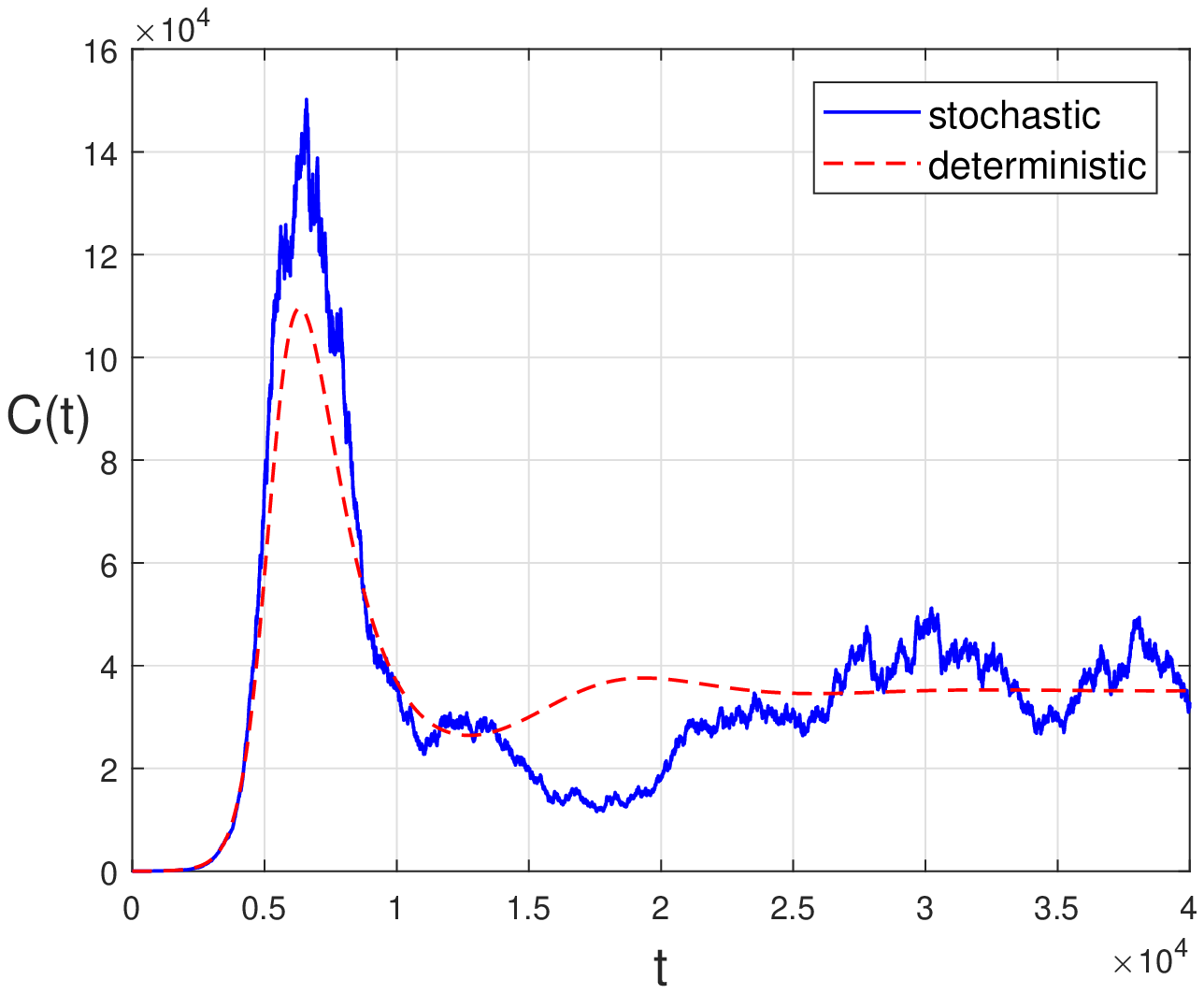}    \includegraphics[width=2in]{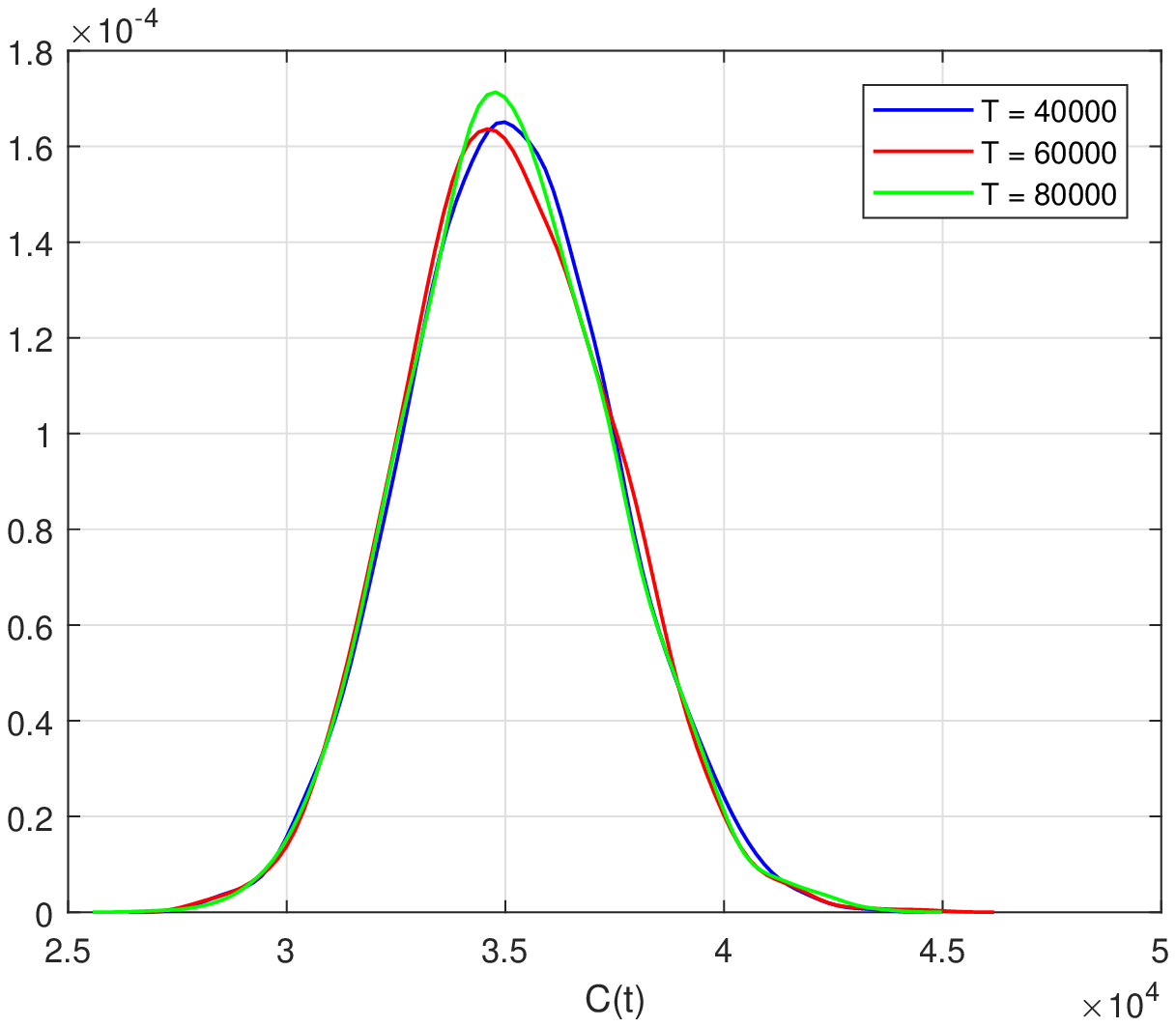}
    \includegraphics[width=2in]{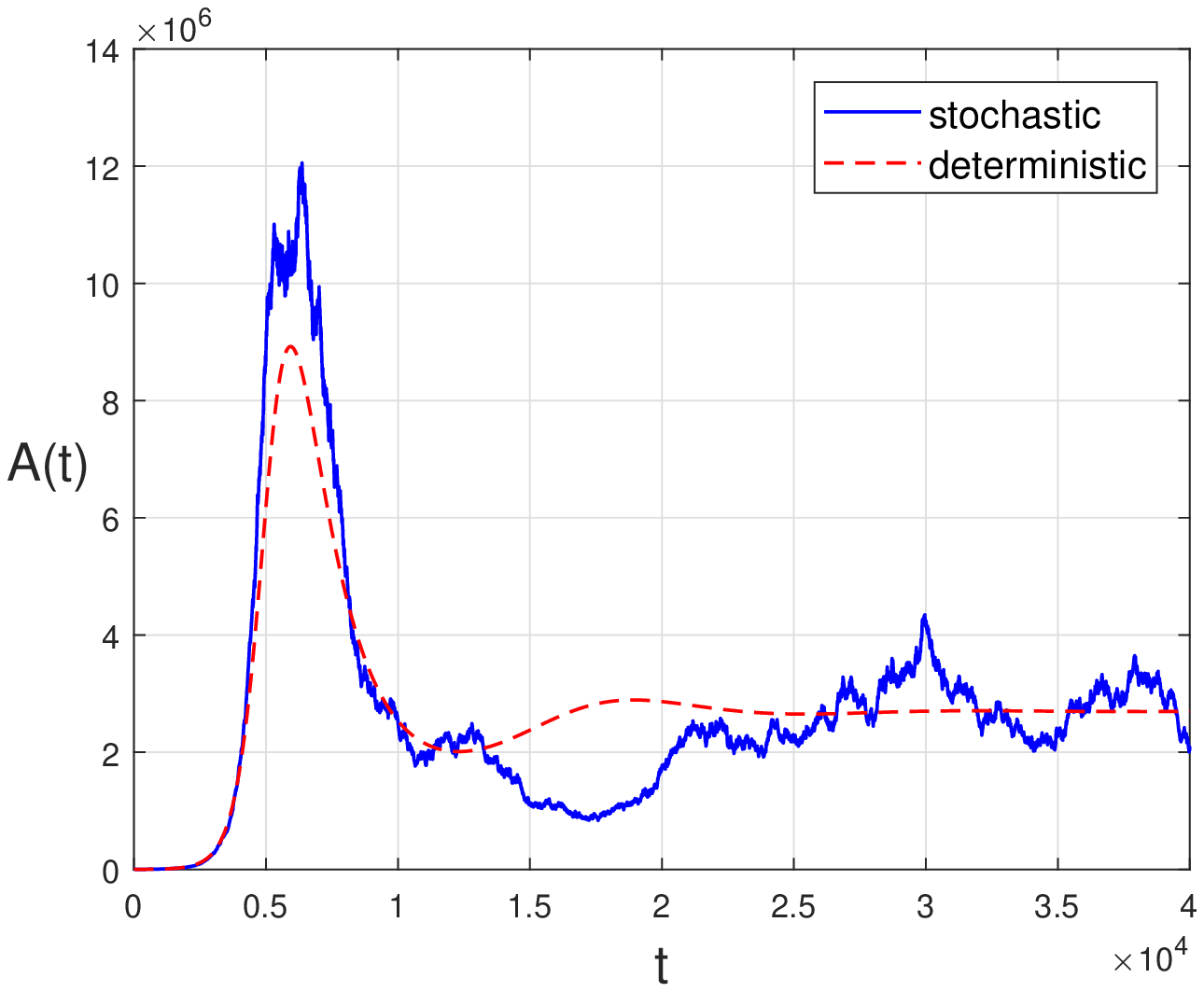}    \includegraphics[width=2in]{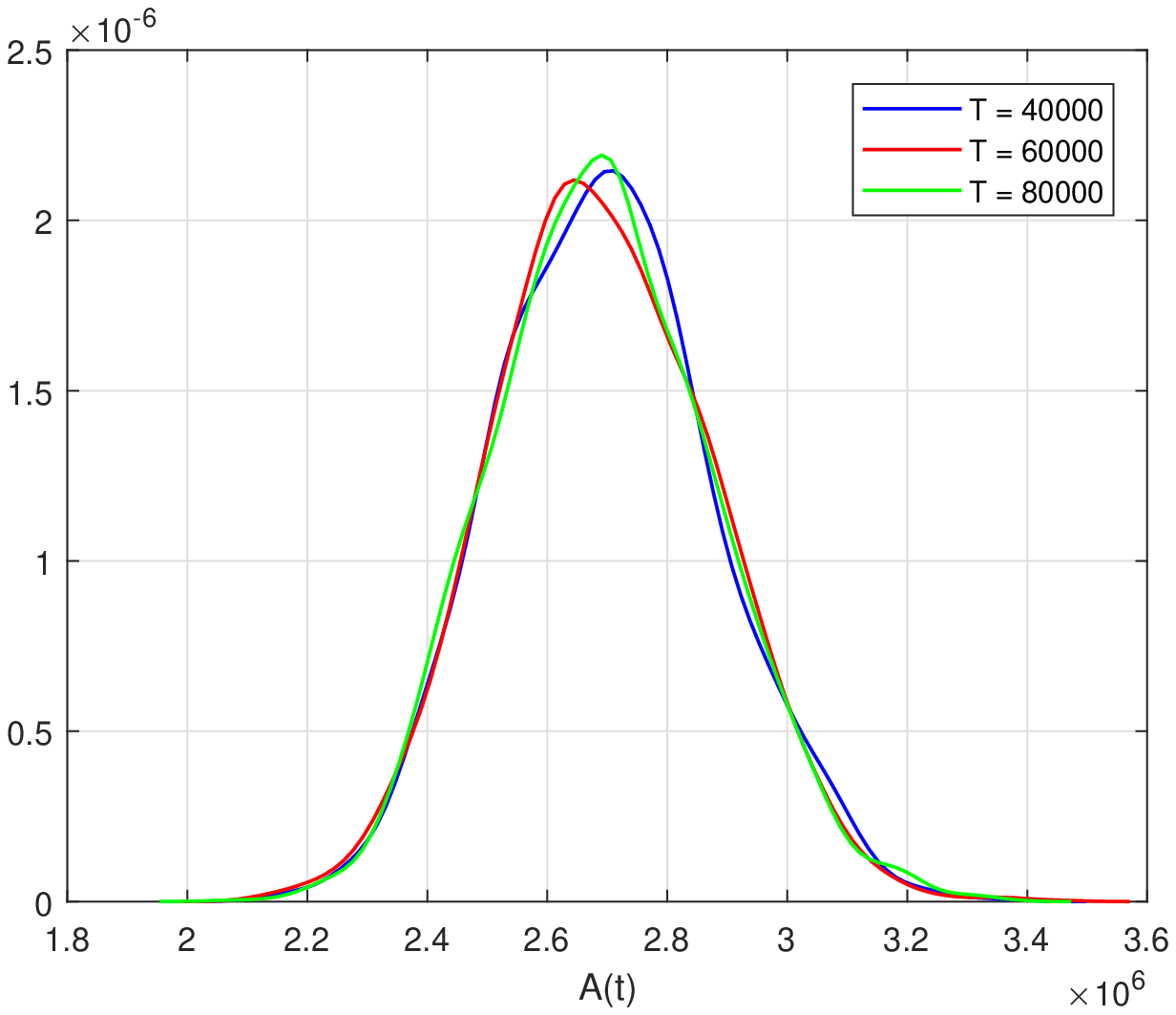}
    \caption*{Figure 6.1\quad The persistence and stationary distribution of $S_u, S_a, I, C, A$ in model (2.3)}
\end{figure}
We set $\beta=\frac{0.1065}{229800000}$ and keep the remaining
parameters and initial values same with those in Figure 6.1. So,
$\hat R_0^e \approx 0.0229<1$ and
\begin{equation}\nonumber
     0.0148 \approx \mu>0.5(\sigma_1^2\vee \sigma_2^2\vee\sigma_3^2\vee\sigma_4^2\vee \sigma_5^2)=0.00125.
\end{equation}
By Theorem 5.1, HIV/AIDS is extinct in Figure 6.2.
\begin{figure}[H]
    \centering
    \includegraphics[width=2in]{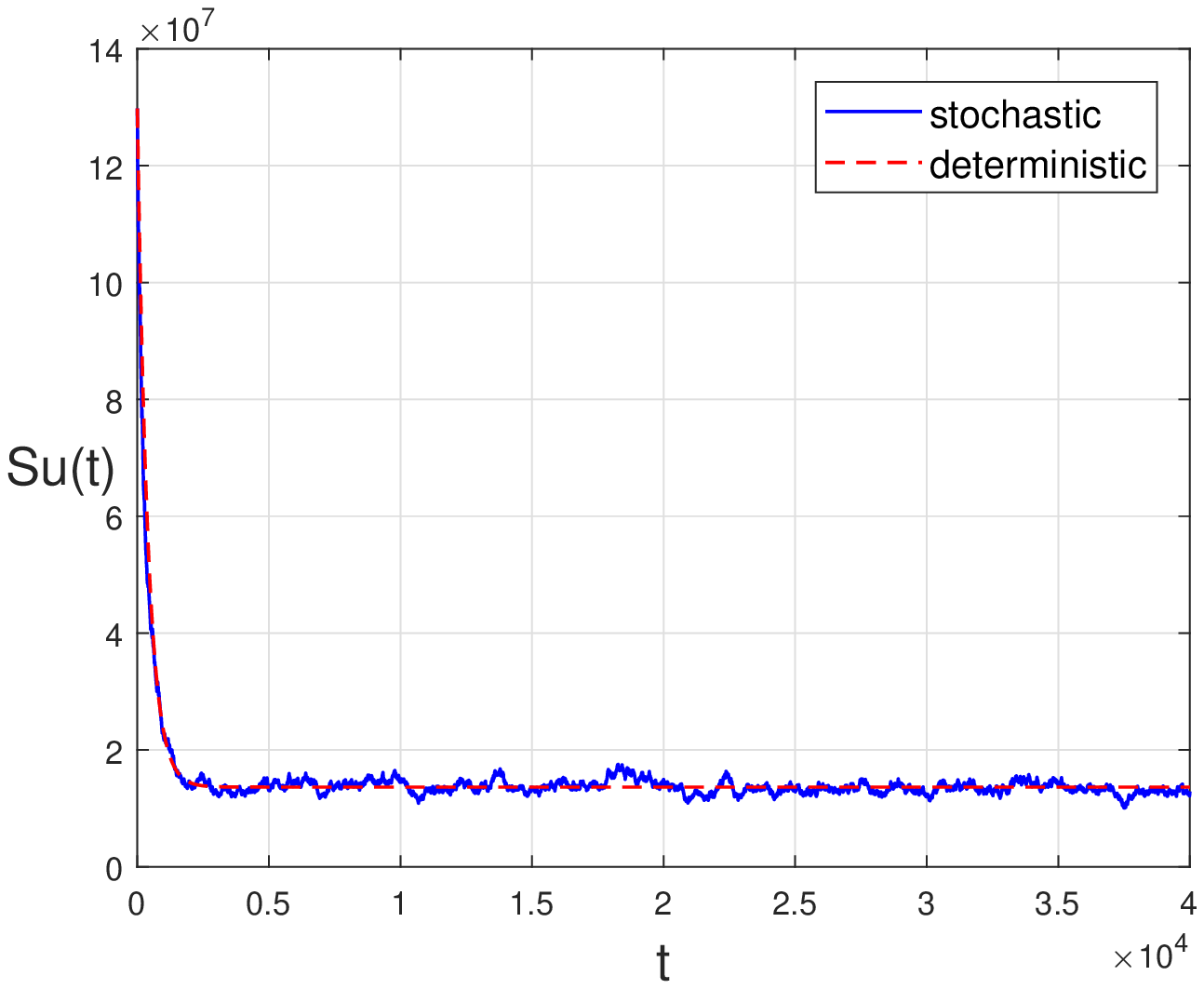}
    \includegraphics[width=2in]{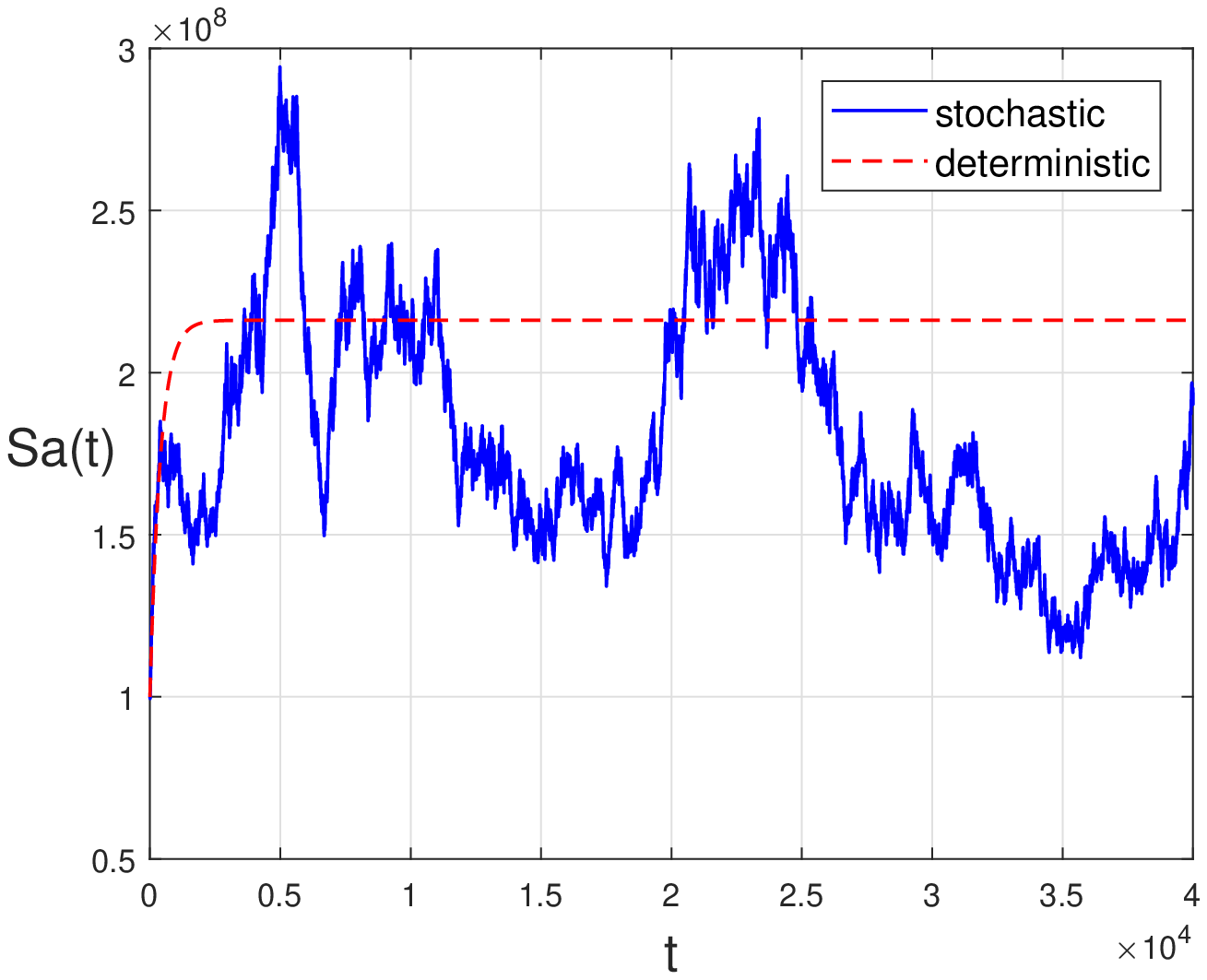}
    \includegraphics[width=2in]{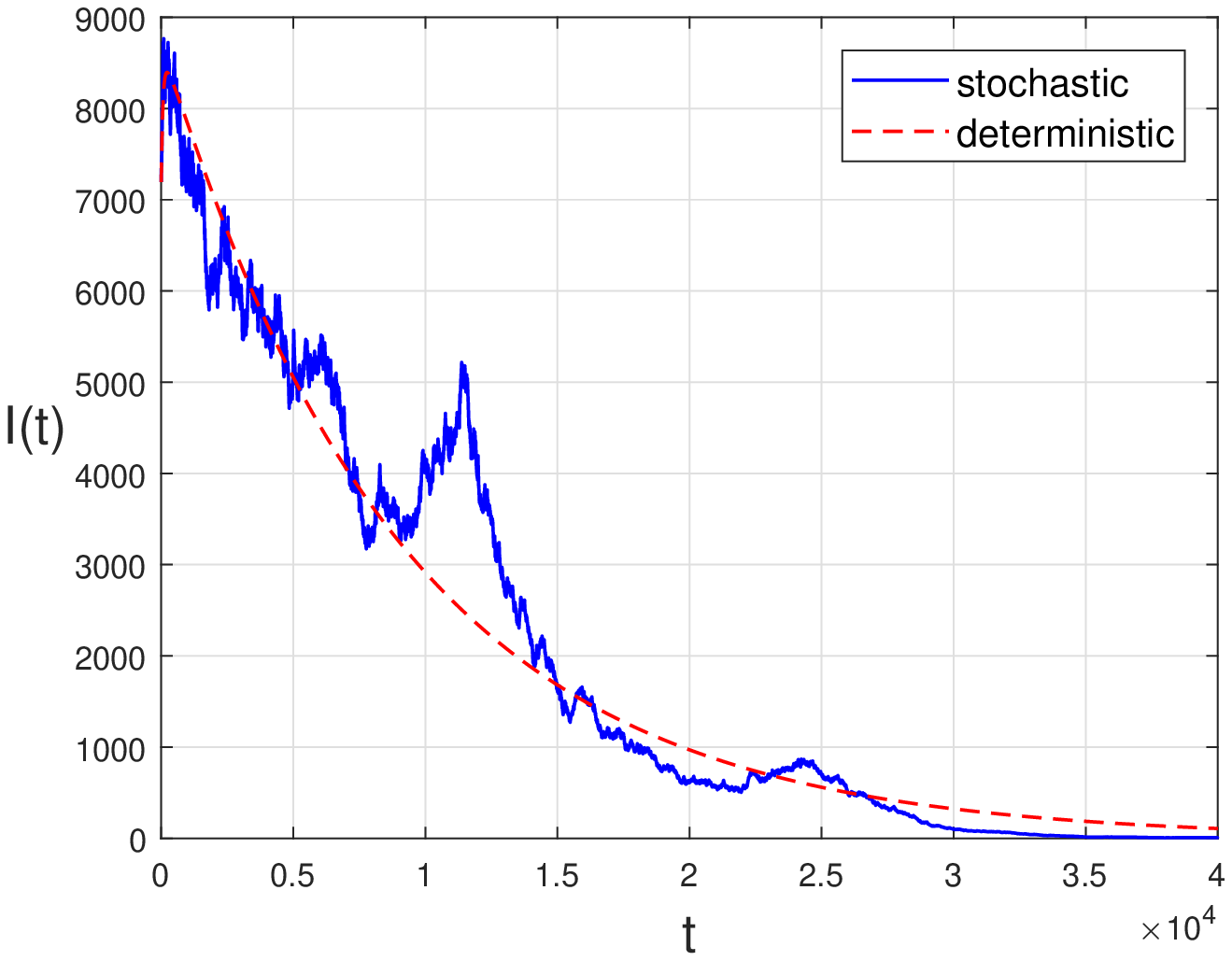}
    \includegraphics[width=2in]{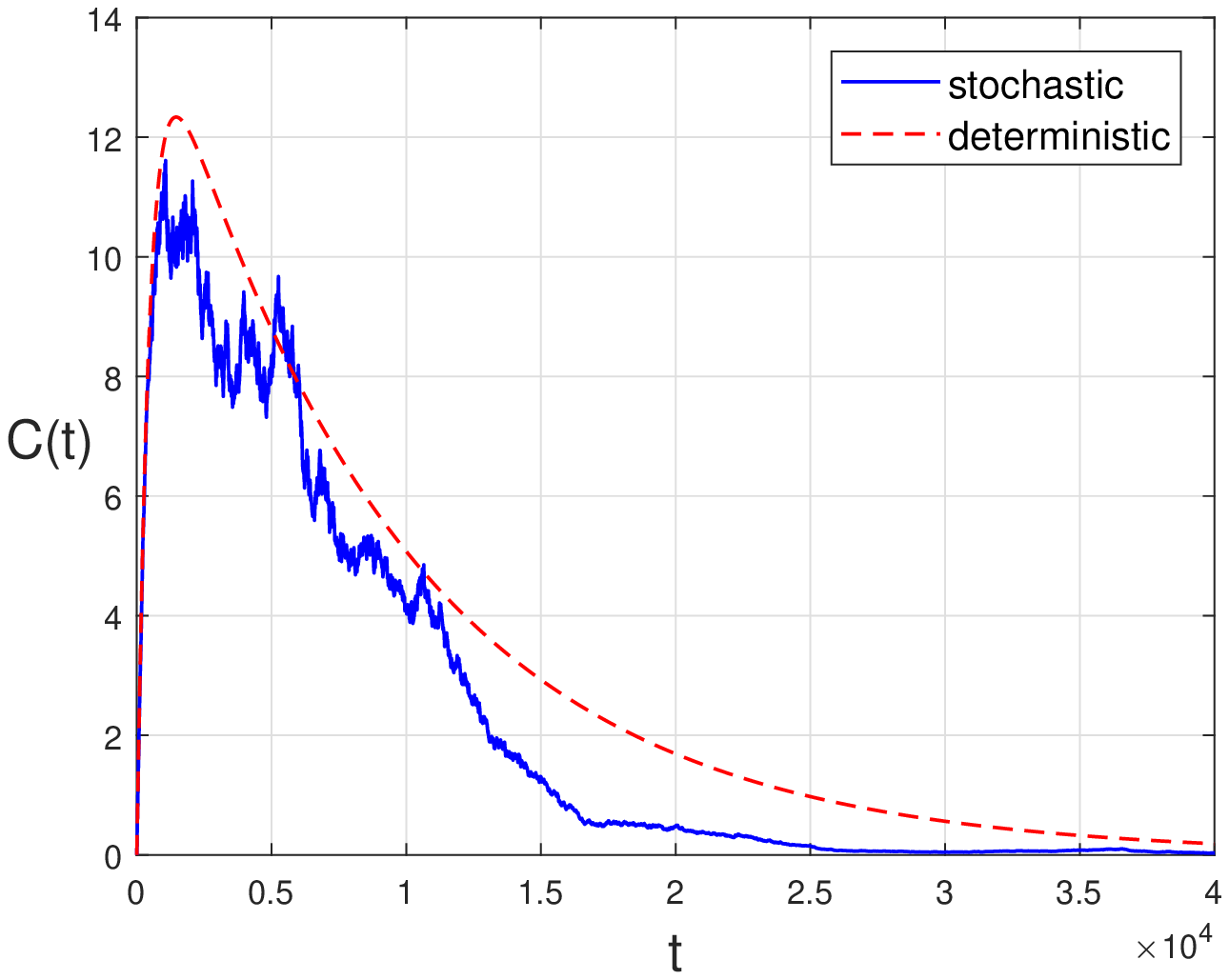}
    \includegraphics[width=2in]{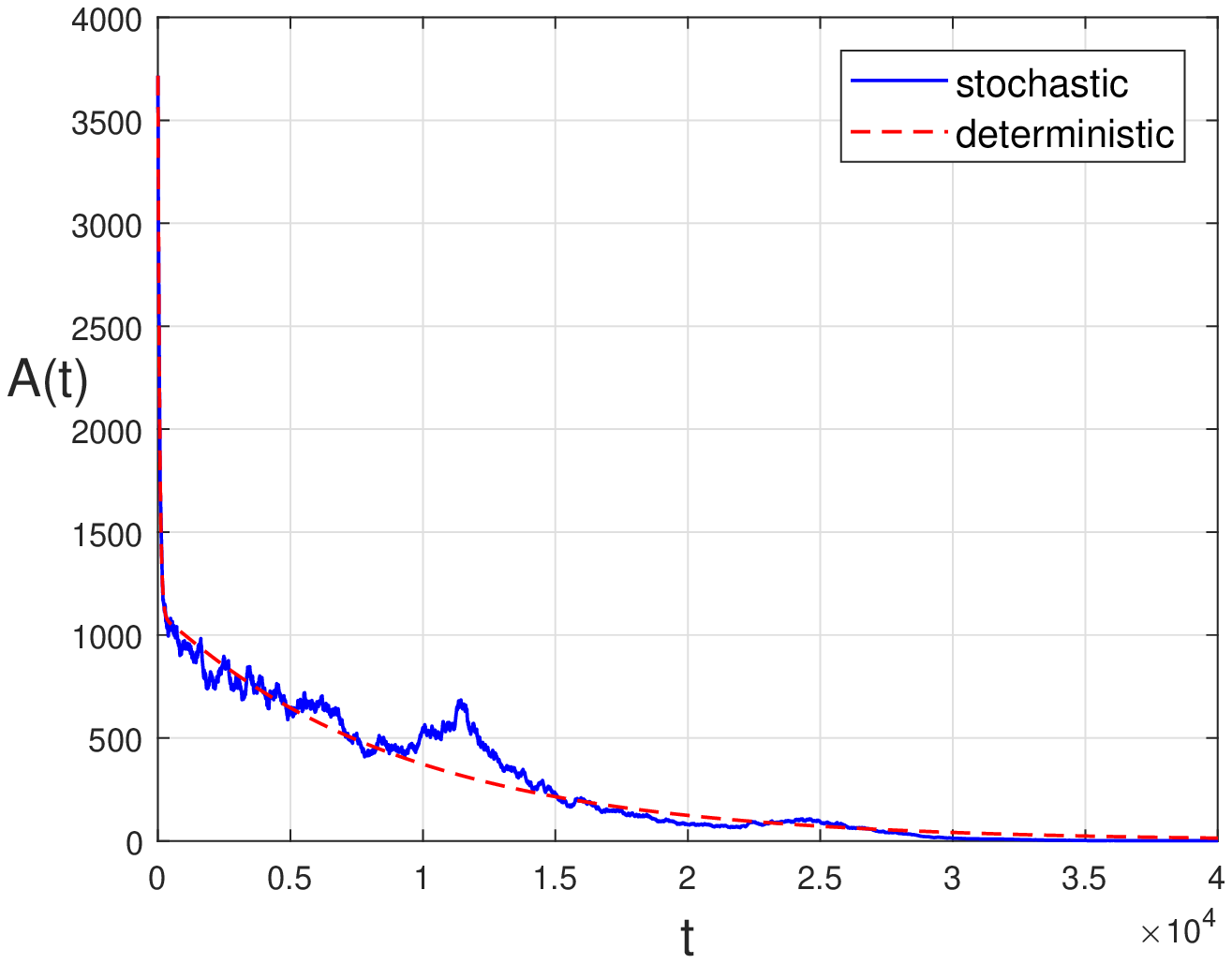}
    \caption*{Figure 6.2\quad The extinction of HIV/AIDS}
\end{figure}

Next, we discuss the impacts of the parameter $\varepsilon$, and
let $\sigma_i(i=1, 2, 3, 4, 5)\\=0.02$ and other parameters and the
initial values be same with those in Figure 6.1. We observe that
HIV/AIDS is persistent as $\varepsilon$ increases, while the
population size of the infected decreases significantly in Figure
6.2 Therefore, the enhancement of $\varepsilon$ is of significant
importance for the prevention and control of HIV/AIDS.
\begin{figure}[H]
    \centering
    \includegraphics[width=2in]{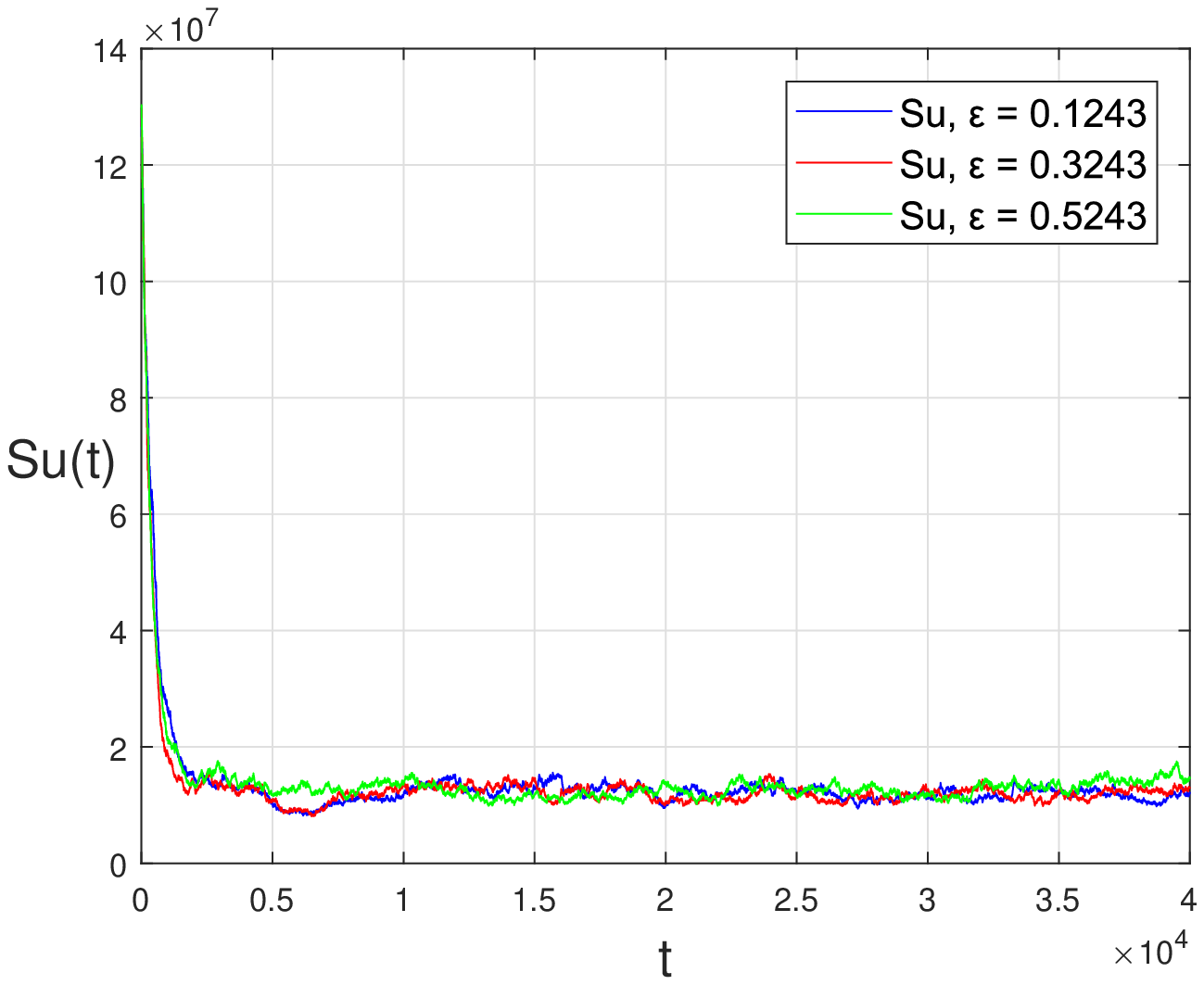}
    \includegraphics[width=2in]{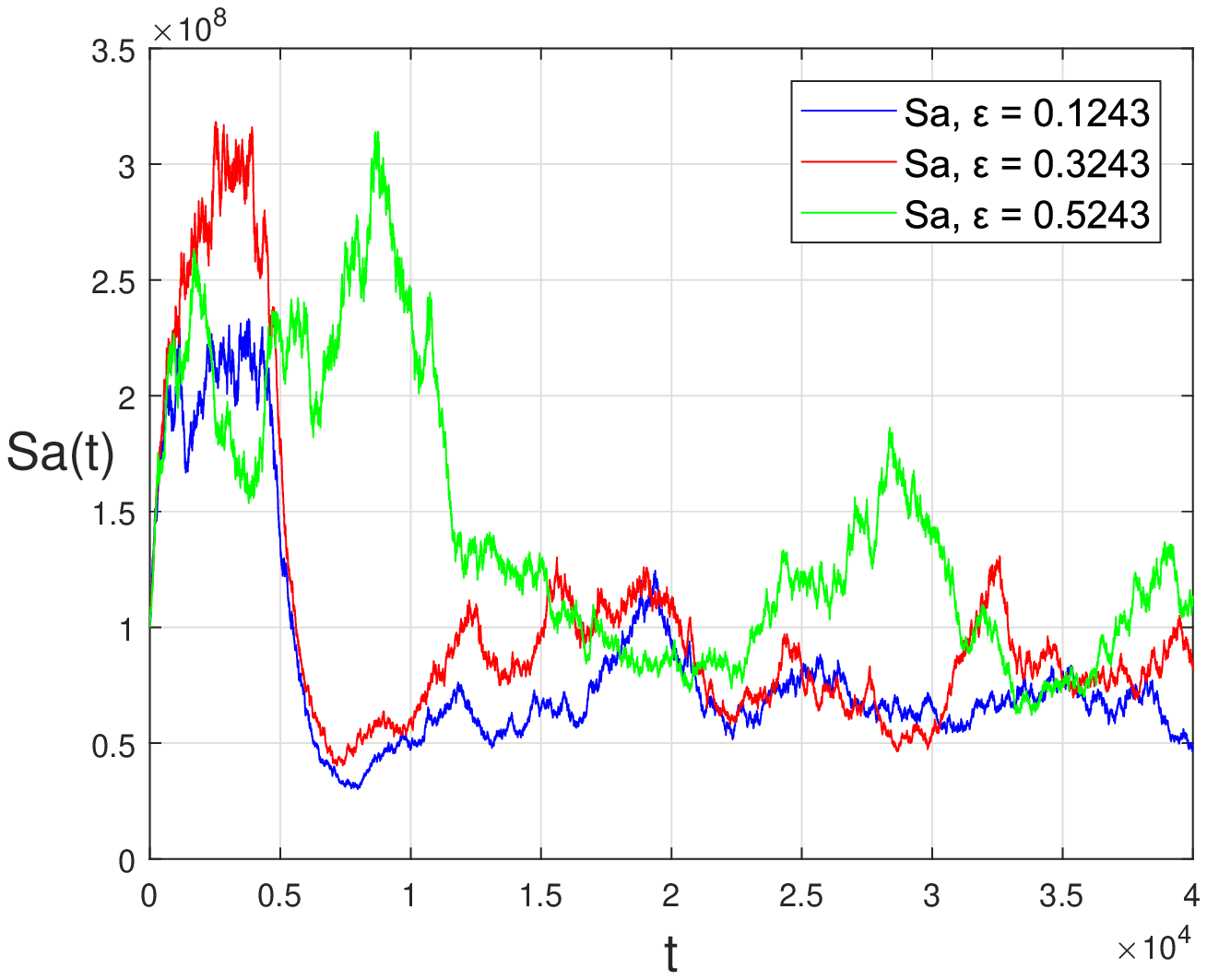}
    \includegraphics[width=2in]{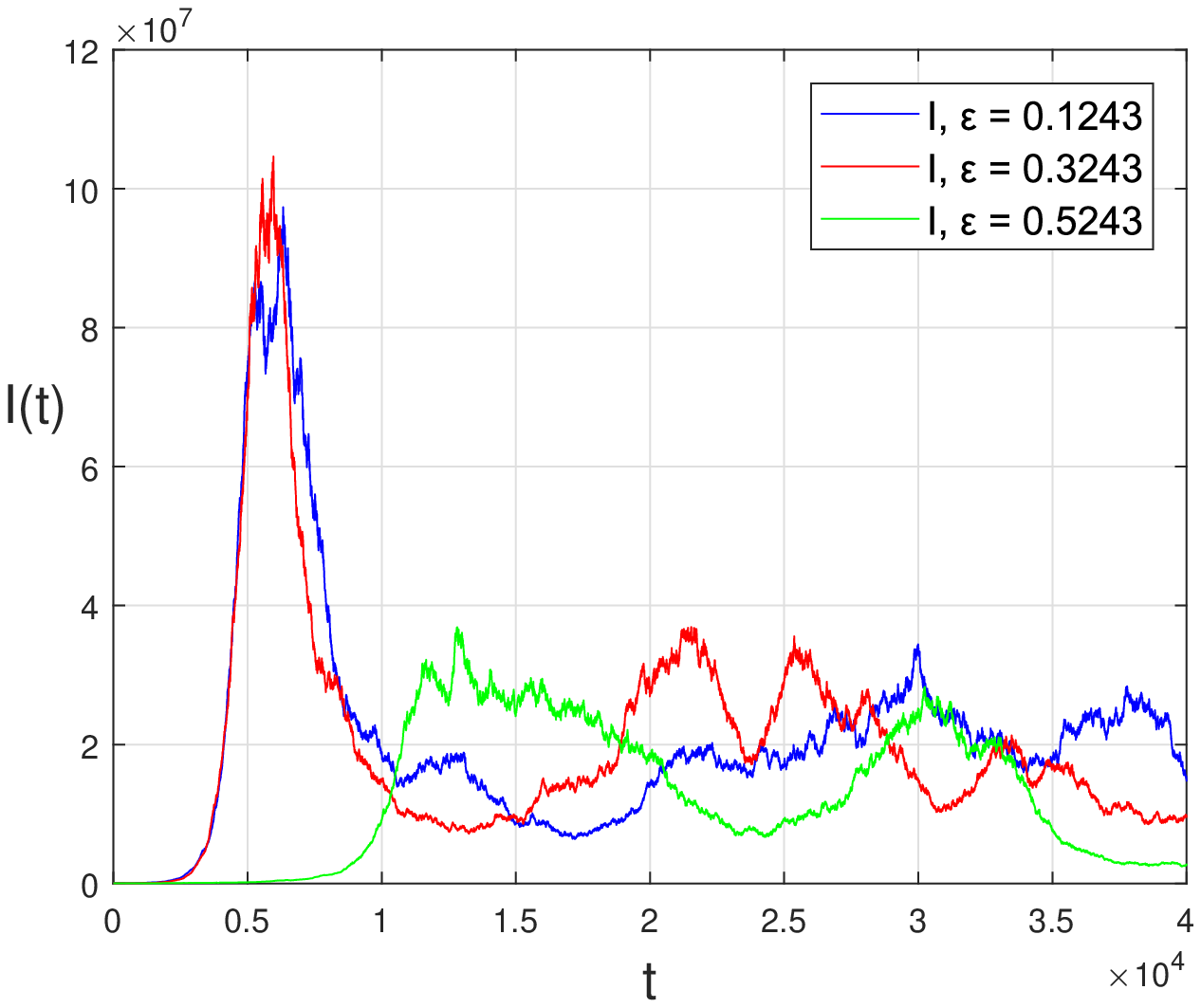}
    \includegraphics[width=2in]{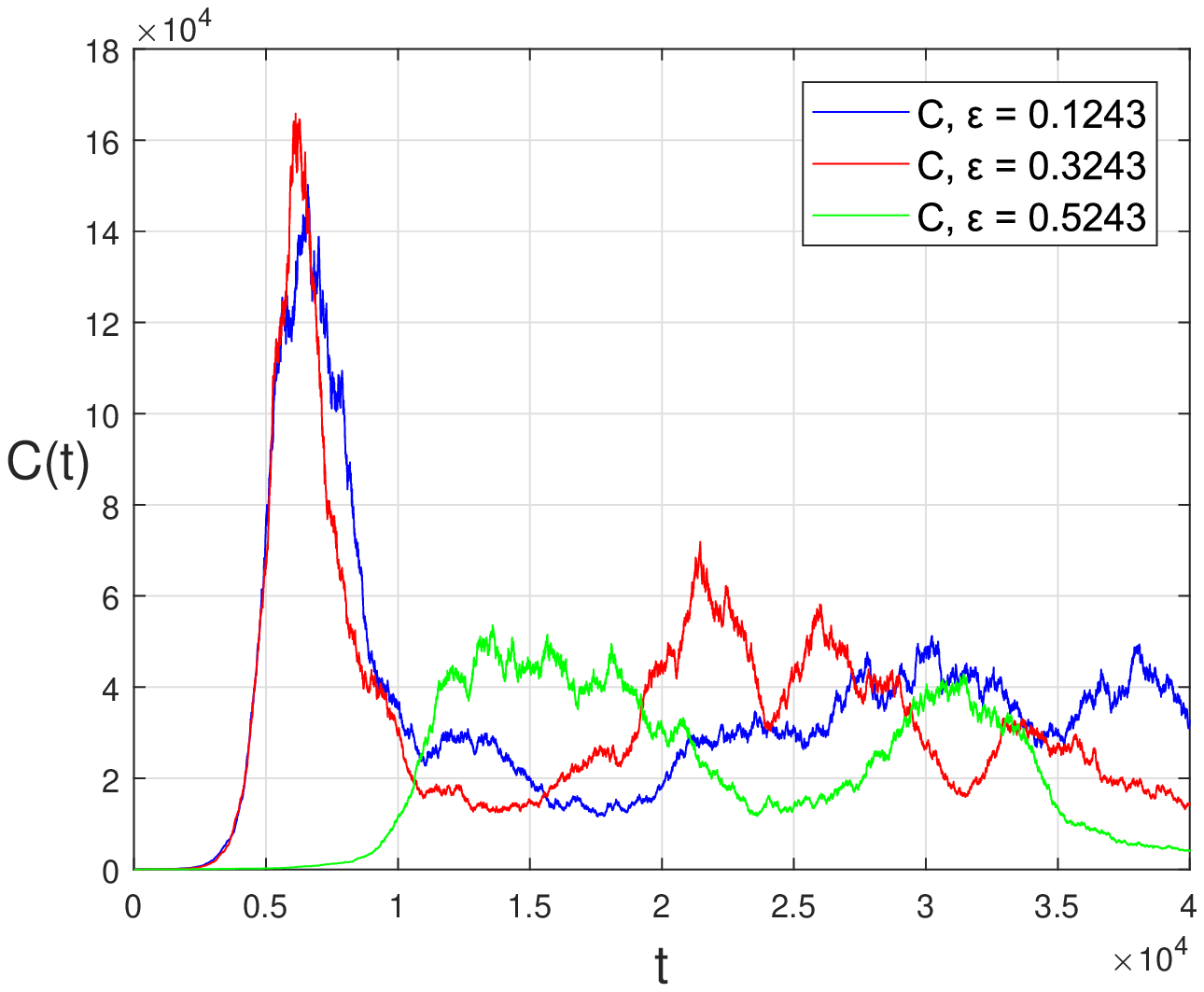}
    \includegraphics[width=2in]{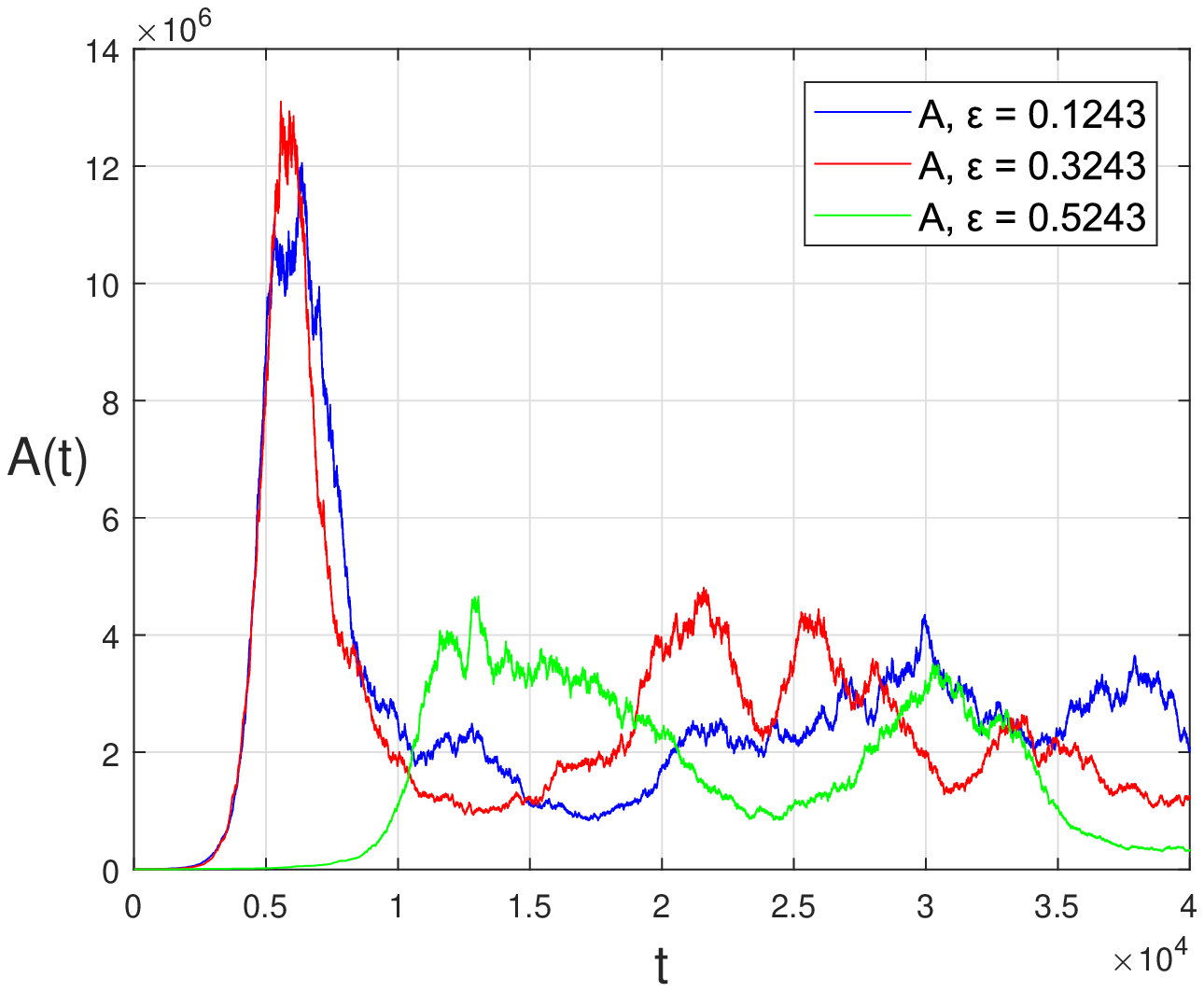}
    \caption*{Figure 6.3\quad  The impact of $\varepsilon$}
\end{figure}

\noindent\textbf{Example 6.2.}

We perform the numerical
simulations on the spread of HIV/AIDS in China for next five
decades, and provide some suggestions for the epidemics in this
example. Since the population size for the year 2014 was
1376460000 in China, and the average life span of the population
was 76.34 \cite{ref32}, we assume that the natural growth rate,
the natural mortality rate and the infection rate for the
population are respectively
$$
\lambda=\frac{1376460000}{76.34},\mu=\frac{1}{76.34},
\beta=\frac{0.71}{1376460000}.
$$
By the data in 2014 in Zhao et al.\cite{ref29}, the number of the
individuals with HIV/AIDS was 500579, and the number of the
individuals with ART was 295358, therefore we assume that
$$S_u=1088230000, S_a=288230000, I=153193, C=295358, A=52128$$
and choose $\Delta = 10^{-2}$ and $\alpha=0.13,  \varepsilon=0.5,
\eta=0.18, \upsilon=0.72, \gamma=0.14, \rho=0.82, \delta=0.42 $ as
other parameters for simulation.

Firstly, we collect the data for the individuals with HIV/AIDS
from the year 2014 to 2020 in China in Zhao et al.\cite{ref29}
(year 2014-2018), Liu et al.\cite{ref30} (year 2019) and He
\cite{ref31} (year 2020). We thus adopt Runge Kutta method to fit
the parameters, and the simulation with fitted parameters and the
data are shown on the left of Figure 6.4. By Theorem 4.1, we
derive that $R_ 0^s=2.8942>1$ as $\sigma_i=0.05(i=1,2,3,4,5)$. Further,
we govern the data in 2020 as the initial values to preform the
simulations for next 5 decades in China as presented on the right
of Figure 6.4. It is easy to observe that although the spread of
HIV/AIDS in China is running in a low epidemic level, there still
exists a risk of the exponential growth for HIV/AIDS control.
\begin{figure}[H]
    \centering
    \includegraphics[width=2in]{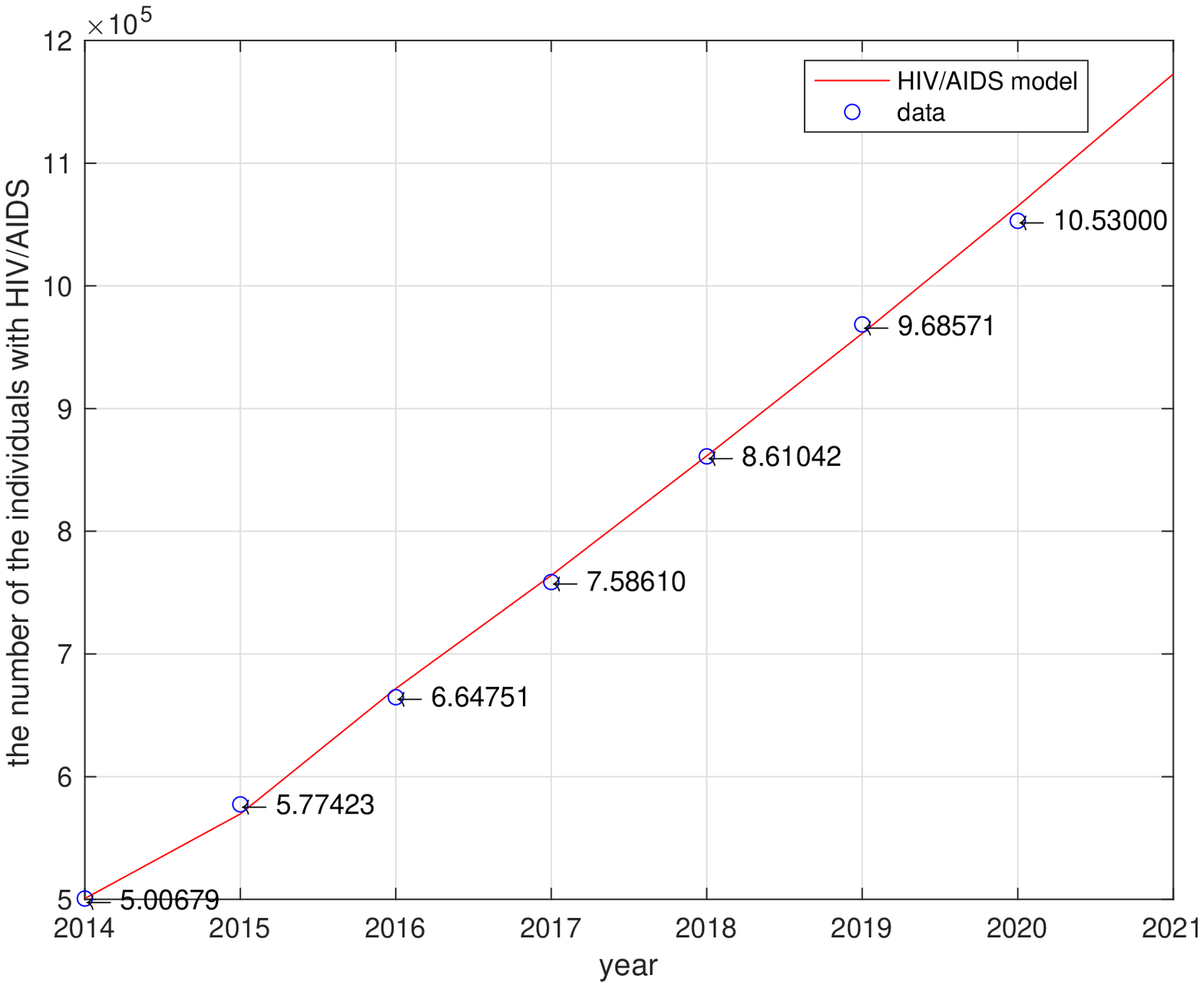}
    \includegraphics[width=2in]{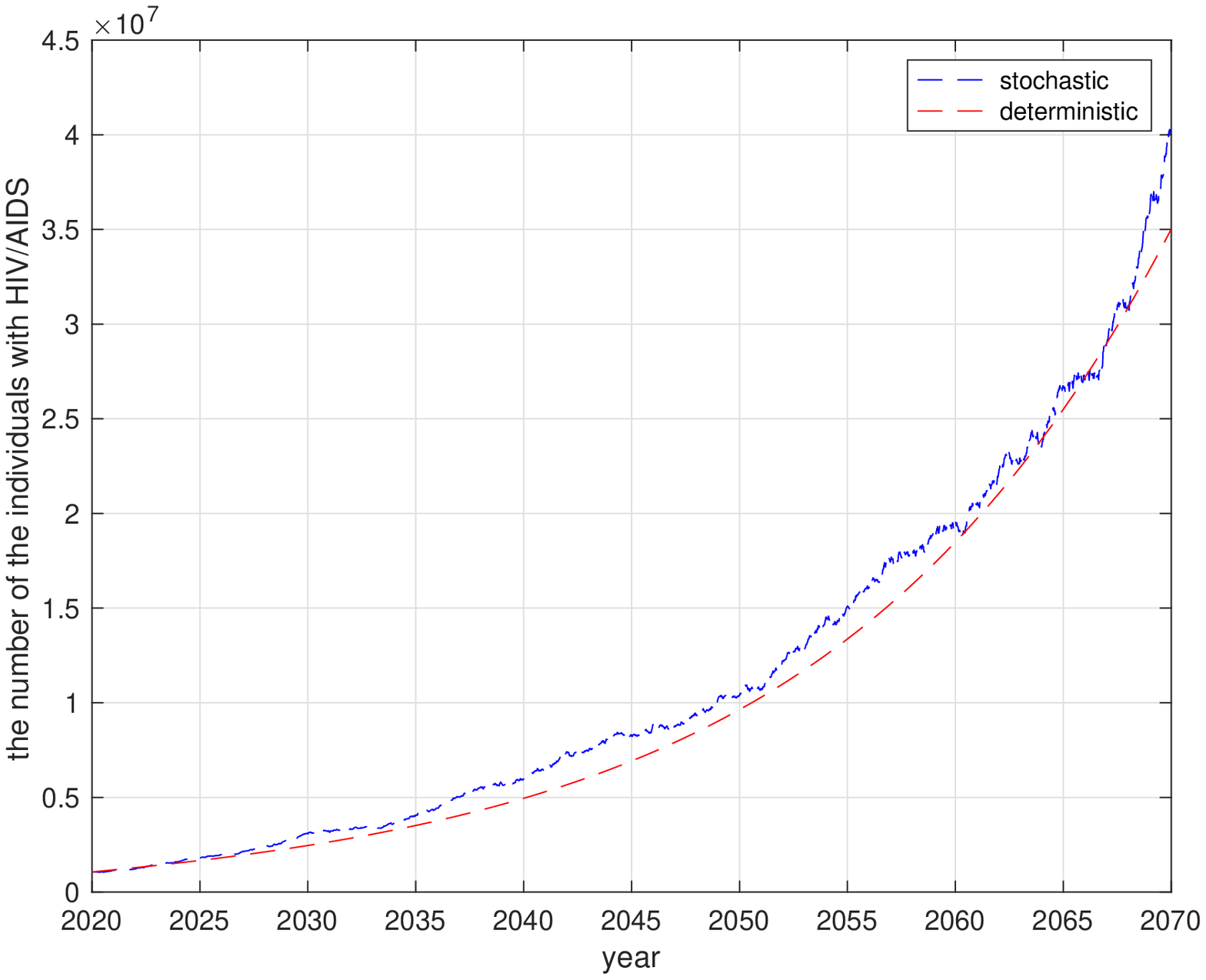}
    \caption*{Figure 6.4\quad  Data fitting for the individuals with HIV/AIDS from 2014 to 2022,
        and prediction of HIV/AIDS for next 5 decades in China}
\end{figure}

Next, we discuss the impacts of main parameters to the control of
HIV/ AIDS in China. The extensive publicity and detailed publicity
on HIV/AIDS are two important ways to prevent and control the
spread of HIV/AIDS in China. In practice, the extensive publicity
improves the number of the individuals having protection awareness
from $S_u$ to $S_a$ by varying $\alpha$, which further reveals
that less impact on the transmission of HIV/AIDS occurs as
$\alpha$ increases (see the left in Figure 6.5). Meanwhile, the
detailed publicity presents more details for the individuals who
are infected by HIV, and the isolations within 72 hours are
usually adopted to reduce the infection rates, which further
suppresses the number of the individuals with HIV/AIDS as
$\varepsilon$ increases (see the right in Figure 6.5).

Therefore, the extensive publicity and the detailed publicity
including lessons and lectures of AIDS in universities and
communities to the target population play significant roles to
prevent the spread of HIV/AIDS.
\begin{figure}[H]
    \centering
    \includegraphics[width=2in]{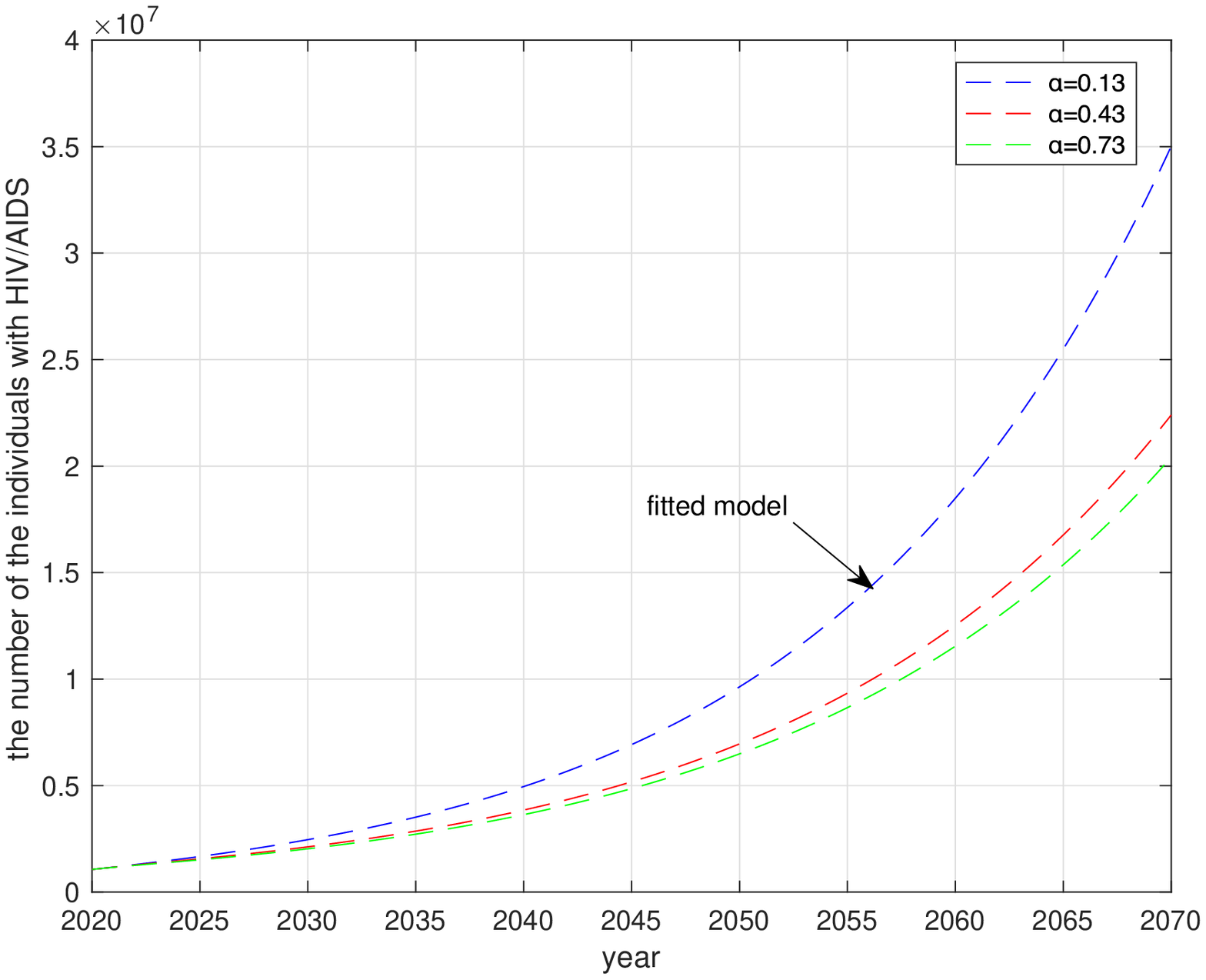}
    \includegraphics[width=2in]{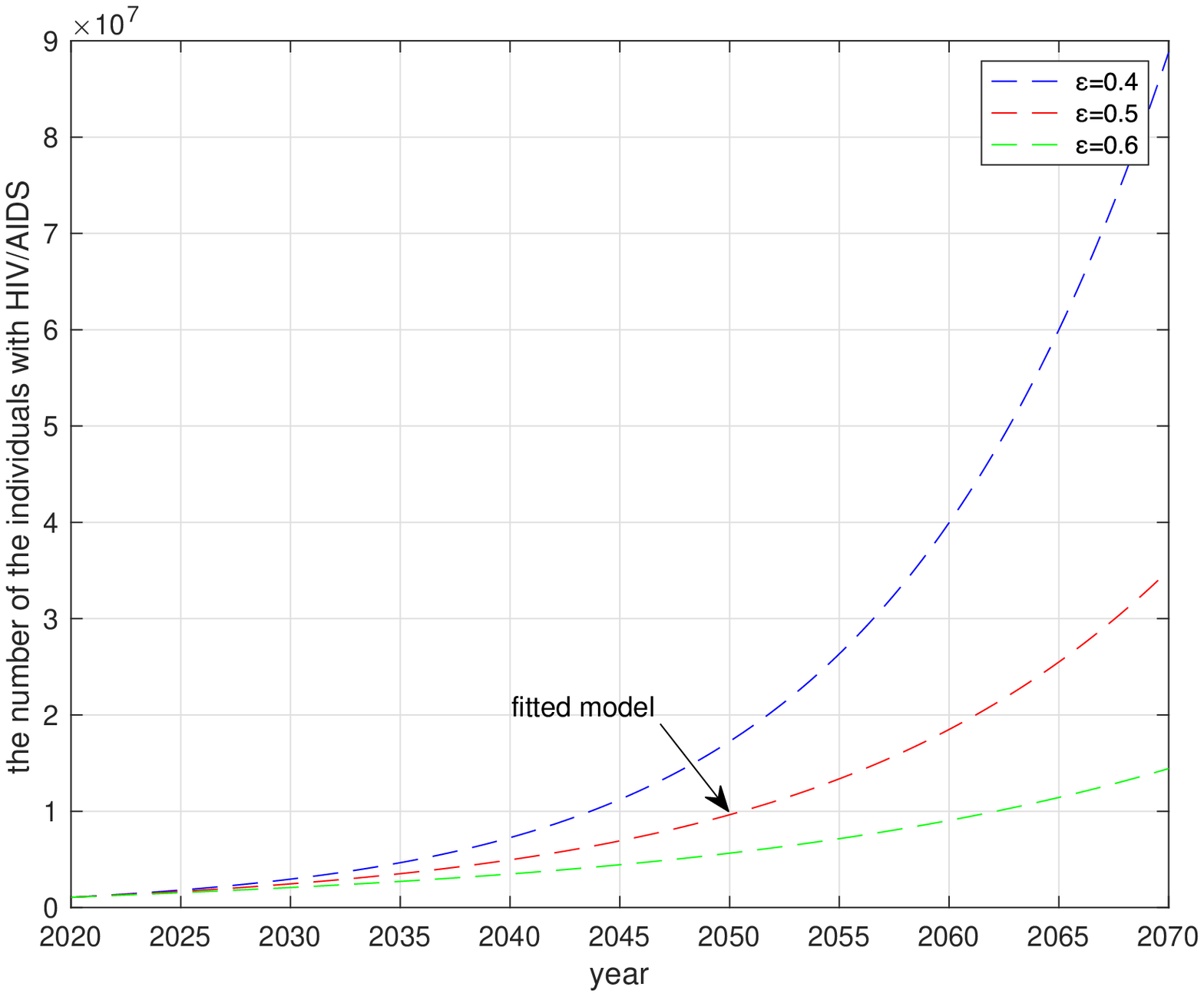}
    \caption*{Figure 6.5\quad Impacts of the extensive publicity and detailed publicity on HIV/AIDS}
\end{figure}

The prompt and continuous antiretroviral therapy (ART) after being
infected is helpful to each individual with HIV/AIDS. Figure 6.6
demonstrates the simulations with distinct values of $\eta$, which
also verify that ART suppresses the rapid growth of the
individuals who are with HIV/AIDS as $\eta=0.08$.
\begin{figure}[H]
    \centering
    \includegraphics[width=2in]{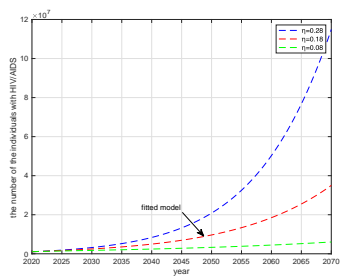}
    \caption*{Figure 6.6\quad Impact of the continuous ART }
\end{figure}

We also notice that the transmission rates $\gamma$ and $v$ affect
the long-term epidemics of HIV/AIDS in China. More precisely, when
$\gamma$ increases (also the period that the individuals with
HIV/AIDS see the doctors in hospital and get checked becomes
shorter), the number of the individuals with HIV/AIDS decrease.
Meanwhile, when $v$ increases (also the period that the
individuals with AIDS stay at hospital becomes shorter), the number
of the individuals with HIV increases as presented in Figure 6.7.
\begin{figure}[H]
    \centering
    \includegraphics[width=2in]{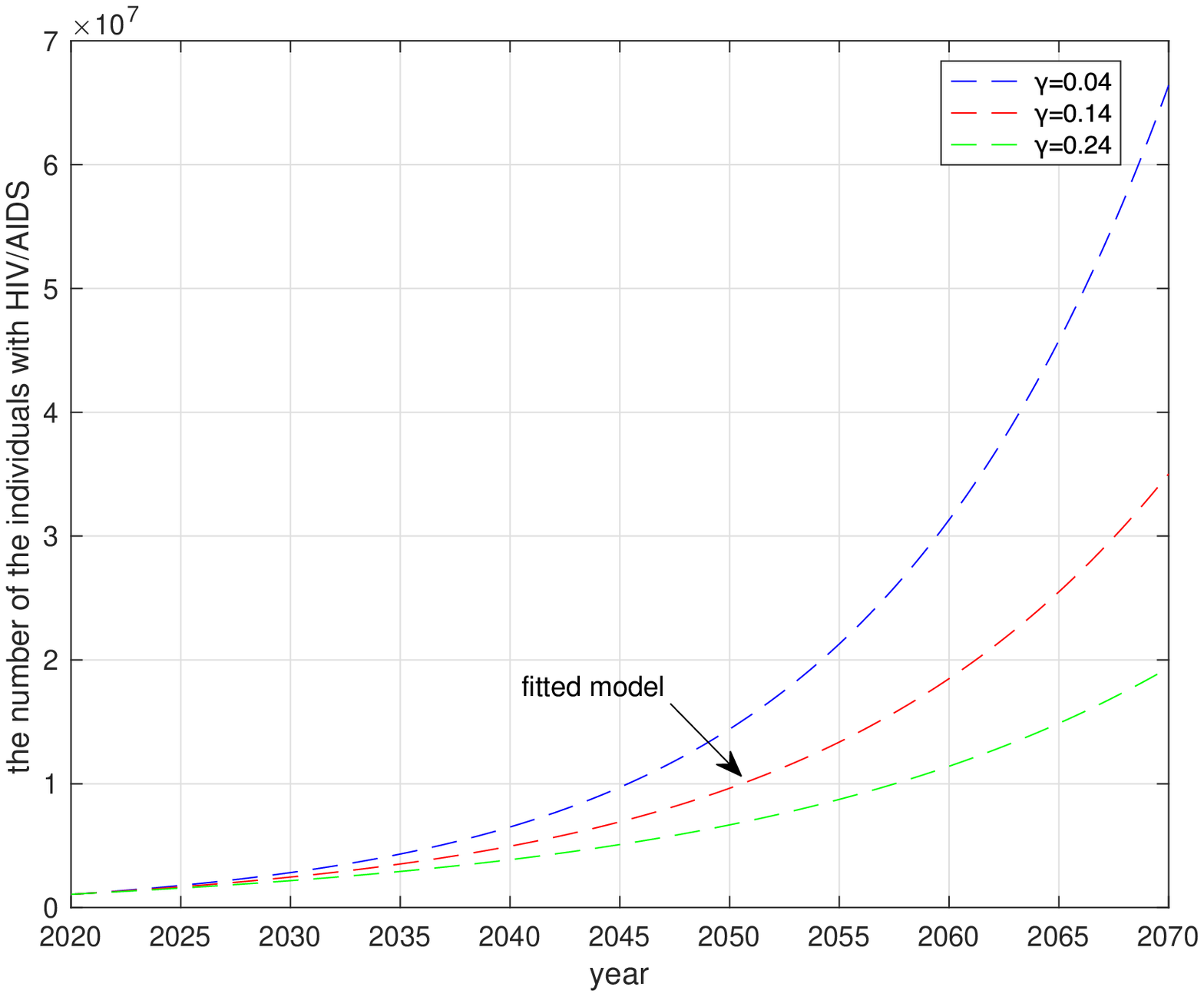}
    \includegraphics[width=2in]{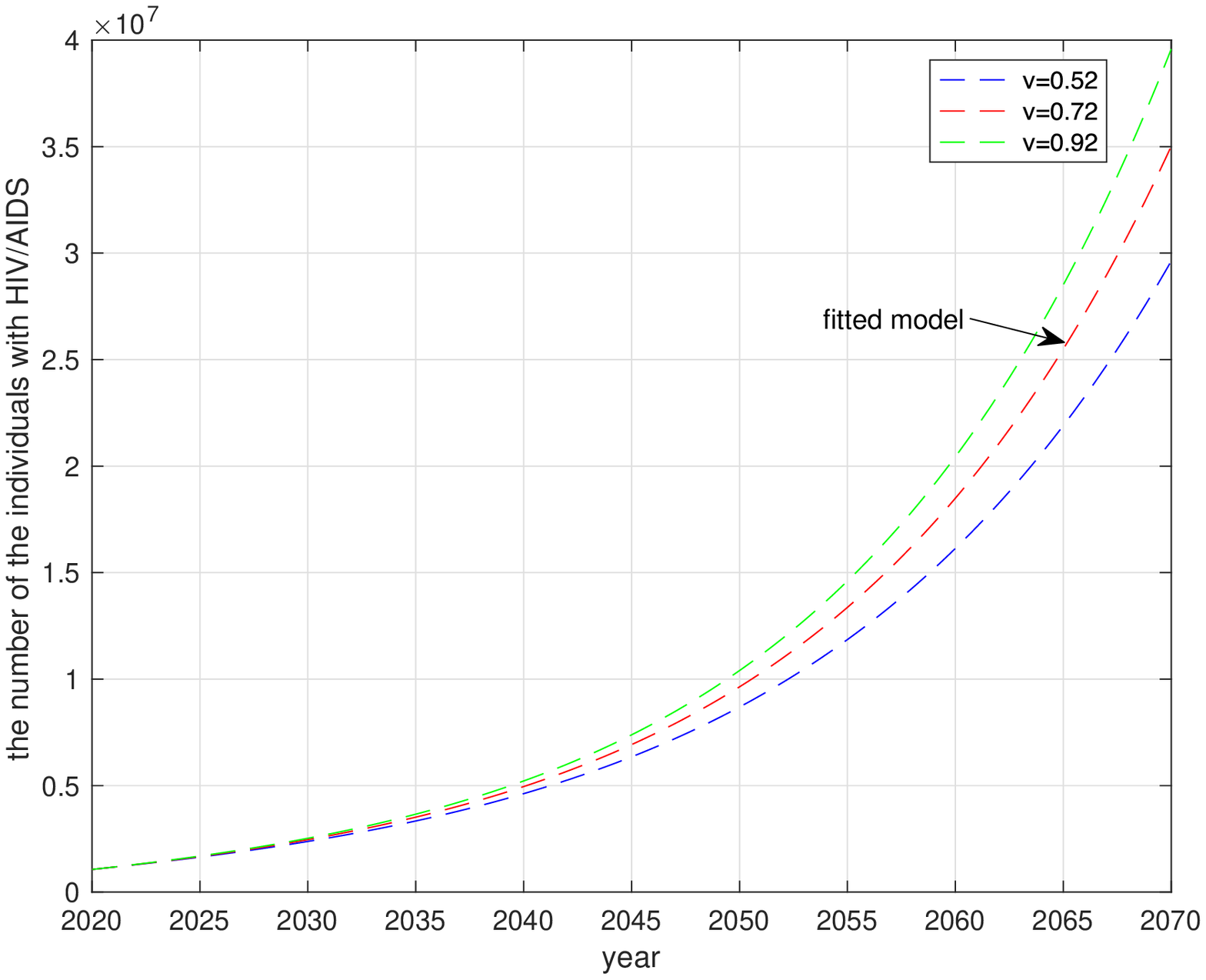}
    \caption*{Figure 6.7\quad  Impacts of $\gamma$ and $v$ on HIV/AIDS}
\end{figure}

\section{Conclusions and discussions}
We propose a stochastic epidemic model with the bilinear incidence
rate and show that the existence of a global positive solution. By
constructing Lyapunov functions, we also show that the stochastic
model has an ergodic stationary distribution when $R^s_0 > 1$.
Moreover, the sufficient condition $R^e_0 < 1$ for the extinction
of HIV/AIDS is obtained. The corresponding simulations verify that
the numbers of the individuals in Indonesia and in China decreas
when the detailed publicity and the continuous ART are governed to
prevent and control the spread of HIV/AIDS. Therefore, we suggest
that all countries should enhance the systematic and detailed publicity
on AIDS, which are of significant importance to control the growth
of HIV/AIDS. For instance, taking the isolations within 72 hours and
receiving prompt ART treatment after being infected are the
effective measurements for the elimination of AIDS by 2030.

Compared with the conclusions derived by Fatmawati et al.
\cite{ref5}, the disease-free equilibrium point $P_0$ attracts the
solution of (2.3) under condition $R_0^e<1$ (see Theorem 5.1),
which also means that HIV/AIDS becomes extinct as the intensities
of the white noises increase. Meanwhile, the solution of (2.3)
fluctuates around the endemic equilibrium $P^*$, and the solution
has a unique ergodic stationary distribution when $R_0^s>1$ (see
Theorem 4.1).

We also point out that the expressions for $R_0^s$ and $R_0^e$ are
two distinct indices for indicating the prevalence of HIV/AIDS.
When  the intensities of the white noises disappear, $R_0^s$ turns
into $R_0$ in (2.2). Theoretically, model (2.3) has a smaller
index for the persistence of HIV/AIDS  than that of model (2.1).
The long-term properties of model (2.3) are quite different from
the results in Fatmawati et al. \cite{ref5}.

\begin{acknowledgements}
\noindent\vskip10pt The research of F.Wei is supported in part by
Natural Science Foundation of Fujian Province of China
(2021J01621) and Technology Development Fund for Central Guide
(2021L3018); the research of X.Mao is supported by the
Royal Society, UK (WM160014, Royal Society Wolfson Research Merit
Award), the Royal Society and the Newton Fund, UK (NA160317, Royal
Society-Newton Advanced Fellowship), the EPSRC, The Engineering
and Physical Sciences Research Council ( EP/K503174/1 ) for their
financial support.

\end{acknowledgements}


\section*{Conflict of interest}
All authors consent to publish the main results of this paper on
Journal of Dynamics and Differential Equations. All authors declared
that we did not and do not have any conflicts of interest with any other
institutions and groups.



\end{document}